\documentclass[doublespacing]{elsart}

\usepackage{amssymb}
\usepackage{epsfig}
\usepackage{tabularx}
\usepackage[english]{babel}
\usepackage{graphicx}
\usepackage[numbers]{natbib}
\usepackage{color}
\usepackage[table,xcdraw]{xcolor}
\usepackage{multirow}
\usepackage{epstopdf}
\usepackage{subfigure}
\usepackage{bm}
\usepackage{algorithm}
\usepackage{algpseudocode}
\usepackage{color, soul} %color,  soul package

\usepackage{amsmath}% http://ctan.org/pkg/amsmath
\begin{document}

\setulcolor{blue} %set underlining color
\setstcolor{red} %set overstriking color
\sethlcolor{yellow} %set highlighting color

\begin{frontmatter}

\title{High-resolution transport of regional level sets for evolving complex interface networks}


\author{Shucheng Pan,}  \author{Xiangyu Hu,} \and  \author{Nikolaus A. Adams}
\address[label1]{Lehrstuhl f\"{u}r Aerodynamik und Str\"{o}mungsmechanik, Technische Universit\"{a}t
M\"{u}nchen,
 85748 Garching, Germany}
\begin{abstract}
In this paper we describe a high-resolution transport formulation of the regional level-set approach 
for an improved prediction of the evolution of complex interface networks.
The novelty of this method is twofold:
(i) construction of local level sets and reconstruction of a global regional level sets,
(ii) locally transporting the interface network by employing high-order spatial discretization schemes 
for improved representation of complex topologies.
Various numerical test cases of multi-region flow problems, 
including triple-point advection, single vortex flow, mean curvature flow,
normal driven flow and dry foam dynamics, 
show that the method is accurate and suitable for a wide range of complex interface-network evolutions.
Its overall computational cost is comparable to the Semi-Lagrangian regional level-set method while the prediction accuracy is significantly improved. The approach thus offers a \textbf{viable} alternative to previous interface-network level-set method.
\end{abstract}
\begin{keyword}
multi-region problem; interface network; level-set method; interface capturing; high-order scheme
\end{keyword}
\end{frontmatter}

\section{\label{sec:intro}Introduction}
Multi-region problems can occur when the motion of more than two immiscible fluids is to be described. In this case the interface network, separating the different fluid regions,
evolves in time due to interactions of the different fluids across interface segments. These interactions often can be described by local fluid properties. Important applications include shock-driven multiphase flows \cite{betney2015computational, thomas2012drive, haines2014three}, astrophysical events \cite{chen2014two, lentz2015three, wongwathanarat2015three}, foam dynamics \cite{weaire2001physics, biance2011topological, saye2013multiscale, kim2014numerical},
multi-cellular tissue dynamics \cite{hilgenfeldt2008physical, rauzi2008nature, manning2010coaction, marinari2012live, osterfield2013three},
and grain coarsening in polycrystalline materials
\cite{geiger2001simulation, krill2002computer, elsey2011large, torres2015numerical}.

A range of numerical models have been proposed
to compute the evolution of interface networks for multi-region problems.
Generally, they can be classified as Lagrangian or Eulerian methods according to the representation of the interface.
With Lagrangian methods, such as front-tracking \cite{unverdi1992front}, immersed-boundary \cite{peskin2002immersed}, or arbitrary Lagrangian-Eulerian (ALE) \cite{hirt1974arbitrary} methods,
the interface is represented explicitly by conforming discretization elements.
Although these methods have been extended to multi-region systems
\cite{brakke1992surface, galera2010two, loubere2010reale, kucharik2010comparative,
kim2010numerical, kim2014numerical},
it is difficult to handle complex topological changes
during interface-network evolution, especially in three dimensions.
With Eulerian methods, such as volume-of-fluid (VOF) \cite{hirt1981volume}
and level-set methods \cite{osher1988fronts},
the interface is reconstructed from scalar fields, i.e., volume fraction or level-set field.
The interface is represented implicitly and captured by solving the corresponding transport equations.
These methods generally can handle complex interface evolution with topological changes
and are straightforward to implement in three dimensions.
However, they often exhibit low accuracy due to numerical dissipation introduced by the transport-equation discretization.

Two additional problems are encountered when Eulerian methods are applied to multi-region problems \cite{merriman1994motion, zhao1996variational, starinshak2014new}. One is that the number of scalar fields increases with the number of regions
and entails additional computational operations and memory cost \cite{chan2007some}. The other is that the interface reconstruction can produce voids and overlaps where more than two regions meet \cite{zhao1996variational, chan2007some, starinshak2014new}. There are two main level-set-based approaches for capturing the evolution of multiple regions. One is to define multiple level-set functions (referred to as multiple level-set method in this paper) and to assign these to corresponding regions \cite{merriman1994motion, zhao1996variational, starinshak2014new}, followed by solving separate level-set transport equations for each region. Different numerical procedures may be employed to prevent the generation of voids and overlaps during the interface reconstruction. For example, Starinshak et al. \cite{starinshak2014new} propose $N(N-1)/2$ level-set functions to represent all interfaces of $N$ regions and an additional pairwise voting strategy to remove overlaps and voids. The algorithm copes with the interface reconstruction problem but requires that a larger number of different level-set fields are stored and evolved in time than with multiple level-set methods \cite{merriman1994motion, zhao1996variational}. The method of Vese and Chan \cite{vese2002multiphase} reduces the number of level-set fields from $N$ to $\log_2 N$ \citep{chan2007some} and it avoids the generation of voids and overlaps. Naturally, using a single level-set function to represent an arbitrary number of regions is the optimal strategy to address the memory overhead \citep{chan2007some}. Lie et al. \cite{lie2006variant} and Chung and Vese \cite{chung2009image} develop such methods for image segmentation. Also, the regional level-set method \cite{zheng2009simulation} addresses the problem of multiple level-set fields \cite{kim2010multi} by employing a combination of a single unsigned level-set field and an integer region indicator function. Interface reconstruction for more than two regions, is handled by employing a low-order Semi-Lagrangian scheme. We note that this low-order Semi-Lagrangian scheme is more dissipative than high-order finite-difference schemes. As further development of the regional level-set method, the Voronoi implicit interface method \cite{saye2011voronoi} uses a transported $\epsilon$-level-set and reconstructs the interface network by a reinitialization step based on Voronoi diagrams.

Our objective is to develop an improved regional level-set method that has computational cost comparable to Semi-Lagrangian regional level-set methods \cite{zheng2009simulation, kim2010multi}, but significantly improves prediction accuracy. Specifically, the method inherits the advantages of the original level-set method and the Semi-Lagrangian regional level-set method but is modified in such a way that it is suitable for high-resolution finite-difference discretization of the level-set transport equations. We define local level-set fields to capture the evolution of the interface network by a simple construction operator, followed by a reconstruction operator to obtain the global regional level-set field. The multi-region system and the definition of the regional level-set field are revisited in Sec. \ref{sec2}. The proposed multi-region method is detailed in Sec. \ref{sec3}. Accuracy and robustness are assessed in Sec. \ref{sec4}, followed by a brief conclusion in Sec. \ref{sec-conclusion}.

\section{Multi-region system and regional level set} \label{sec2}

First we introduce the representation of a multi-region system by implicit functions. Let the domain $\Omega$ be an open set in $\mathbb{R}^d$, let $\mathbf{x}\in \Omega$ be an interior point of $\Omega$, and let $\partial \Omega$ be the boundary of the domain, where $d$ is the spatial dimension. Assuming that there are $\mathcal{N}$ regions within this domain, $\Omega$ is the union of a family of disjoint subsets, $\Omega =  \underset{\chi \in \mathrm{X}}{\bigcup} \Omega^{\chi} $, where $\Omega^{\chi}$ is the subdomain of the region $\chi$, and $\mathrm{X} = \{\chi \in \mathbb{N} | 1\leq \chi \leq \mathcal{N}\}$ is the index set for all regions. The entire multi-region system consists of the following elements:

$\bullet$ The region domain $\Omega^{\chi}$ which contains all interior points of the region $\chi$. \\
$\bullet$ The region boundary $\partial \Omega^{\chi}$ which contains the set of boundary points of $\Omega^{\chi}$.\\
$\bullet$ The pairwise interfaces $\Gamma_{\alpha\beta} = \partial \Omega^{\alpha} \cap \partial \Omega^{\beta}$ ($\alpha, \beta \in \mathrm{X}$ and $\alpha \neq \beta $), each as codimension manifold in $\mathbb{R}^d$ that separates two connected regions.\\
$\bullet$ The interface network $\Gamma = \underset{\alpha, \beta \in \mathrm{X}, \alpha \neq \beta}{\bigcup} \Gamma_{\alpha\beta} $ is the union of all pairwise interfaces.\\
$\bullet$ Multiple junctions (high order junctions) $ \mathrm{J} = \underset{\alpha, \beta \in \mathrm{X}, \alpha \neq \beta}{\bigcap} \Gamma_{\alpha\beta}$ are intersections of pairwise interfaces. In most practically relevant applications they are triple points (2D) or triple lines (3D). \\

If there are only two regions ($\mathcal{N}=2$) in the system, the above elements can be implicitly defined by a level-set field (or signed distance function) \cite{osher1988fronts} $\phi : \mathbb{R}^d \rightarrow \mathbb{R}$. The two region domains, $\Omega^{1}$ (or $\Omega^{+}$) and $\Omega^{2}$ (or $\Omega^{-}$), are identified by the sign of the level-set function: $\Omega^{1} = \{\mathbf{x} \in \mathbb{R}^d | \phi(\mathbf{x}) > 0\}$ and $\Omega^{2} = \{\mathbf{x} \in \mathbb{R}^d | \phi(\mathbf{x}) < 0\}$. The region boundaries of $\Omega^{1}$ and $\Omega^{2}$ coincide and are identical to the pairwise interface and interface network, $ \partial \Omega^1 = \partial \Omega^2 = \Gamma_{\alpha\beta} = \Gamma = \{\mathbf{x} \in \mathbb{R}^d | \phi(\mathbf{x}) = 0\}$. There are no multiple junctions in this system. When more than two regions ($\mathcal{N}>2$) exist, as shown in Fig.\ref{stencil}(a), the regional level-set method \cite{zheng2009simulation, saye2011voronoi} can be used to represent the system implicitly. For the regional level-set method, a mapping $\varphi^{\chi} : \Omega \subset \mathbb{R}^d \rightarrow \mathbb{R}\times \mathbb{N}$ is defined as $\varphi^{\chi}(\mathbf{x})=(\varphi(\mathbf{x}), \chi(\mathbf{x}))$, where $\varphi(\mathbf{x})\ge 0$ is the unsigned distance function and $\chi(\mathbf{x})$ is a positive integer region indicator. A subdomain $\Omega^{\alpha}$ is identified by region indicators, $\Omega^{\alpha} = \{\mathbf{x} \in \Omega | \chi(\mathbf{x}) = \alpha\}$ and the interface network is defined as the zero distance contour $\Gamma = \{\mathbf{x} \in \mathbb{R}^d | \varphi(\mathbf{x}) = 0\}$. However, the region boundary, the pairwise interface and multiple junctions can not be identified directly by regional level-set fields. An additional operation, such as the Voronoi method described in \cite{saye2011voronoi} is required for this purpose. In this paper we propose as an alternative a simple construction operator to define these surface contours, see Sec. \ref{sec3}.
\section{Numerical Method} \label{sec3}

\subsection{Global and local index sets of regions} \label{sec3:0}
In this section, we develop the numerical method for the representation of the multi-region system and the evolution of the interface network. For this purpose we consider the multi-region system defined by
$\varphi^{\chi}_{i,j} = (\varphi_{i,j}, {\chi}_{i,j})$ at the centers of finite-volume cells with indices $i$ in $x_1$ and $j$ in $x_2$ coordinate directions on a two-dimensional uniform Cartesian grid, as shown in Fig. \ref{stencil}(a). An extension of the two-dimensional definitions to three dimensions is straightforward. We recall that globally there are $\mathcal{N}$ regions in the considered domain. Although $\mathcal{N}$ may be very large, locally the number of regions which occur in a neighborhood of a cell is limited to a small integer number. Consequently, in this local small subdomain $\omega$, the complexity of the system is significantly reduced, with the maximum possible number of regions $\mathcal{N}_{\omega}$ being
\begin{equation} \label{infinite_cell}
\max{\mathcal{N}_{\omega}} = \min(2^{d},\mathcal{N}).
\end{equation}
Interface networks typically exhibit only very few finite-volume cells which contain multiple junctions. Thus, multiple signed level-set functions defined only locally to capture the evolution of the interface network allow to reduce the computational effort significantly. The required global and local definitions on a two-dimensional uniform Cartesian grid are as follows: \\
$\bullet$ $\mathcal{N}$ global regions are labelled by the global index set $\mathrm{X}$. The center $\mathbf{x}_{i,j}$ of each finite-volume cell $C_{i,j}$ is assigned to one region by an indicator $\chi_{i,j}$ even though $C_{i,j}$ may contain more than one region. Here $i,j$ are the global indices for cells. \\
$\bullet$ A local subdomain is defined as $\omega \subset \Omega$. The $\mathcal{N}_\omega$ regions contained in this local subdomain are labelled by a local index set $\mathrm{X}_\omega = \{r \in \mathbb{N} | 1\leq r \leq \mathcal{N}_\omega\}$. We map the local region indicators into the global index set $\mathrm{X}$, and identify the corresponding global region indicator as $\chi_r$.
A local subdomain $\Omega^{\alpha}_\omega$ of the global region $\Omega^{\alpha}$ is defined by $\Omega^{\alpha}_\omega = \{\mathbf{x} \in \omega | \chi_{r}(\mathbf{x}) = \alpha, r \in \mathrm{X}_\omega\} \subset \Omega^{\alpha}$. Here, $\omega$ can be a finite-volume cell $C_{i,j}$ or a neighborhood, i.e., a set of cells $C_{k,l}$) of $C_{i,j}$, which will be defined in Sec. \ref{sec3:1}. $(k,l)$ denotes the local index pair for cells in such a neighborhood of $C_{i,j}$.

\subsection{Cell neighborhoods and cell types} \label{sec3:1}
For a finite-volume cell $C_{i,j}$ one can define two types of neighborhoods, i.e. cell sets:
a square-shaped near neighborhood
\begin{equation} \label{near-neighborhood}
V_{s}=\{C_{k,l} |i-1<k<i+1, j-1<l<j+1\},
\end{equation}
and a cross-shaped stencil neighborhood
\begin{equation} \label{stencil-neighborhood}
V_{c}=\{C_{k,l} |i-t<k<i+t,l=j\}\cup \{C_{k,l} |k=i,j-t<l<j+t\},
\end{equation}
where $2t+1$ is the required width of the stencil for the level-set transport discretization schemes. The number of separate regions $\mathcal{N}_C$ contained in cell ${C_{i,j}}$ and the number of region indicators in the square-shaped near neighborhood $\mathcal{N}_{s}$ have the relation $1 \leq \mathcal{N}_C \leq \mathcal{N}_{s} \ll \mathcal{N}$ when $\mathcal{N}$ is large. Accordingly, the local index sets $\mathrm{X}_\omega$ of Sec. \ref{sec3:0} specifically for $V_{s}$ and $V_{c}$ are $\mathrm{X}_s = \{r \in \mathbb{N} | 1\leq r \leq \mathcal{N}_s\}$ and $\mathrm{X}_c = \{r \in \mathbb{N} | 1\leq r \leq \mathcal{N}_c\}$, respectively. Note that $V_{c}$, $V_{s}$, $\mathcal{N}_{C}$ and $\mathcal{N}_{s}$ are defined for each cell $C_{i,j}$ and thus depend on the index pair $(i,j)$. We identify $\chi_1=\chi_{i,j}$ as ``primary indicator'', and other region indicators $\chi_r$, $r \in \mathrm{X}_s, r>1$, if they exist, are called ``secondary indicators''.

A cell is denoted as a ``full cell'' if it is not intersected by an interface,
i.e. $\mathcal{N}_C = 1$,
as a ``two-region cut cell'' if it is intersected by an interface segment shared by two regions,
i.e. $\mathcal{N}_C = 2$,
or as a ``complex-region cut cell'' if the interface segment is shared by more than two regions,
i.e. $\mathcal{N}_C > 2$, see Fig.\ref{stencil}(a).
For details on identifying cell types, please refer to Algorithm \ref{algo:type}. Analogously, according to the value of $\mathcal{N}_{s}$, a cell $C_{i,j}$ is categorized as a cell with full type $V_{s}$ ($\mathcal{N}_{s} = 1$), as a cell with two-region type $V_{s}$ ($\mathcal{N}_{s} = 2$), or as a cell with complex-region type $V_{s}$ ($\mathcal{N}_{s} > 2$). We clarify that the square-shaped near neighborhood $V_{s}$ serves as a search stencil to identify all local regions that may share the cell $C_{i,j}$ and to generate the local index set $\mathrm{X}_s$.
On the cross-shaped stencil neighborhood $V_{c}$,
$\mathcal{N}_{s}$ auxiliary local level-set fields (see Sec.\ref{sec3:2})
corresponding to the local regions identified in $V_{s}$
are constructed in order to evolve the interface network by discretized evolution equations
and represent the region boundaries from the zero level set.

\subsection{Description of the local construction and reconstruction operators} \label{sec3:2}
For any cell $C_{k,l}$ in the neighborhood of the cell $C_{i,j}$, we can apply a construction operator for generating the local multiple signed level-set fields and a reconstruction operator for reconstructing the global regional level-set field from the local level-set fields, see Fig. \ref{operator}.
\subsubsection{The local construction operator}
The construction operator is used to generate the $\mathcal{N}_{s}$ local level-set fields $\phi^{r}_{k,l} := \phi_{k,l}(\chi_r)$ from the regional level-set field $\varphi^{\chi}_{k,l}$ at the center of cell $C_{k,l}$ for each region indicator $\chi_r \in \mathrm{X}_s$. The construction operator $\mathbf{C}_r : \mathbb{R}\times \mathbb{N} \rightarrow \mathbb{R}$ is defined as
\begin{equation}\label{construction}
\phi^{r}_{k,l}= \mathbf{C}_r \left(\varphi^{\chi}_{k,l}\right) = \begin{cases}
\varphi_{k,l} \quad {\rm if}~ \chi_{k,l} = \chi_r \cr
-\varphi_{k,l} \quad \rm{otherwise}	
\end{cases},
\end{equation}
and the local level-set function $\phi^{r}(\mathbf{x}) := \phi(\mathbf{x},\chi_r)$ at any point can be obtained by interpolation.

Upon construction of the local multiple signed level-set fields, the normal direction and curvature are obtained at the finite-volume cells,
\begin{equation} \label{normal}
\mathbf{n}= \frac{\nabla \phi^{r}}{|\nabla \phi^{r}|},
\end{equation}
\begin{equation} \label{curvature}
\kappa =\nabla \cdot \frac{\nabla \phi^{r}}{|\nabla \phi^{r}|}.
\end{equation}
Unlike with the original level-set method, the normal direction always points away from the interface network. The local construction operator $ \mathbf{C}_r$ may be inaccurate for curvature calculation when the signed level set is not strictly a distance function, which may happen near multiple junctions, as demonstrated in \ref{sec:2dmean}. This issue can be handled by a construction operator $\mathbf{C}_r^{\star}$ on $V_{s}$ which additionally invokes a re-initialization procedure for recovering the signed-distance property of the local multiple signed level-set fields
\begin{equation}\label{construction}
\phi^{r}_{k,l} = \mathbf{C}_r^{\star} \left(\varphi^{\chi}_{k,l}\right) = \begin{cases}
\mathbf{C}_r \left(\varphi^{\chi}_{k,l}\right) \quad {\rm if}~ \chi_{k,l} = \chi_r \cr
\mathbf{C}_r \left( \min\left(\Delta x, d_s(k,l)\right) \right) \quad \rm{otherwise}	
\end{cases},
\end{equation}
where $d_s(k,l)$ is the distance of the cell center $\mathbf{x}_{k,l}$ from the local interface segment in $V_{s}$.

Note that discrete derivatives in the local level-set advection equation Eq. (\ref{local-advection}), see below, which operate on constructed local level-set functions in the stencil neighborhood $V_{c}$ now can be calculated by high-resolution spatial discretization schemes.

\subsubsection{The local reconstruction operator}
The regional level-set field at an arbitrary point can be reconstructed from the multiple auxiliary local level sets for all region indicators in $V_{s}$ by the reconstruction operator $\mathbf{R} : \mathbb{R}^{\mathcal{N}_{s}} \rightarrow \mathbb{R}\times \mathbb{N}$ given by
\begin{equation}\label{reconstruction}
\varphi^{\chi}(x,y) = \mathbf{R} \left(\phi^{r}(x,y), r \in \mathrm{X}_s\right)
=  \left(\left|\max \phi^{r}(\mathbf{x})\right|, \arg\max_{\chi_r} \phi^{r}(\mathbf{x})\right).
\end{equation}
Note that if the interface network is static, application of the construction operation followed by the reconstruction operation leaves the original regional level-set field invariant. Consider the example in Fig. \ref{operator}. If the regional level set at the central cell (colored by red) is $\varphi^{\chi} = (0.2, 4)$, upon local construction, the $3$ generated local level-set data $\{\phi^{r}|r=1,2,3\}$ at this point are $\{0.2, -0.2, -0.2\}$. Upon applying the reconstruction, we obtain the original regional level-set data $\varphi^{\chi} = (0.2, 4)$.

This reconstructed field is unique. The reconstruction generates topologies without artificial overlaps or voids. Thus it is suitable for determining the regional level set when the interface network has been evolved in time through updating multiple local level-set fields.

\subsection{Evolution of the interface network} \label{sec3:3}
\begin{figure}[p]
\begin{center}
\includegraphics[width=0.95\textwidth]{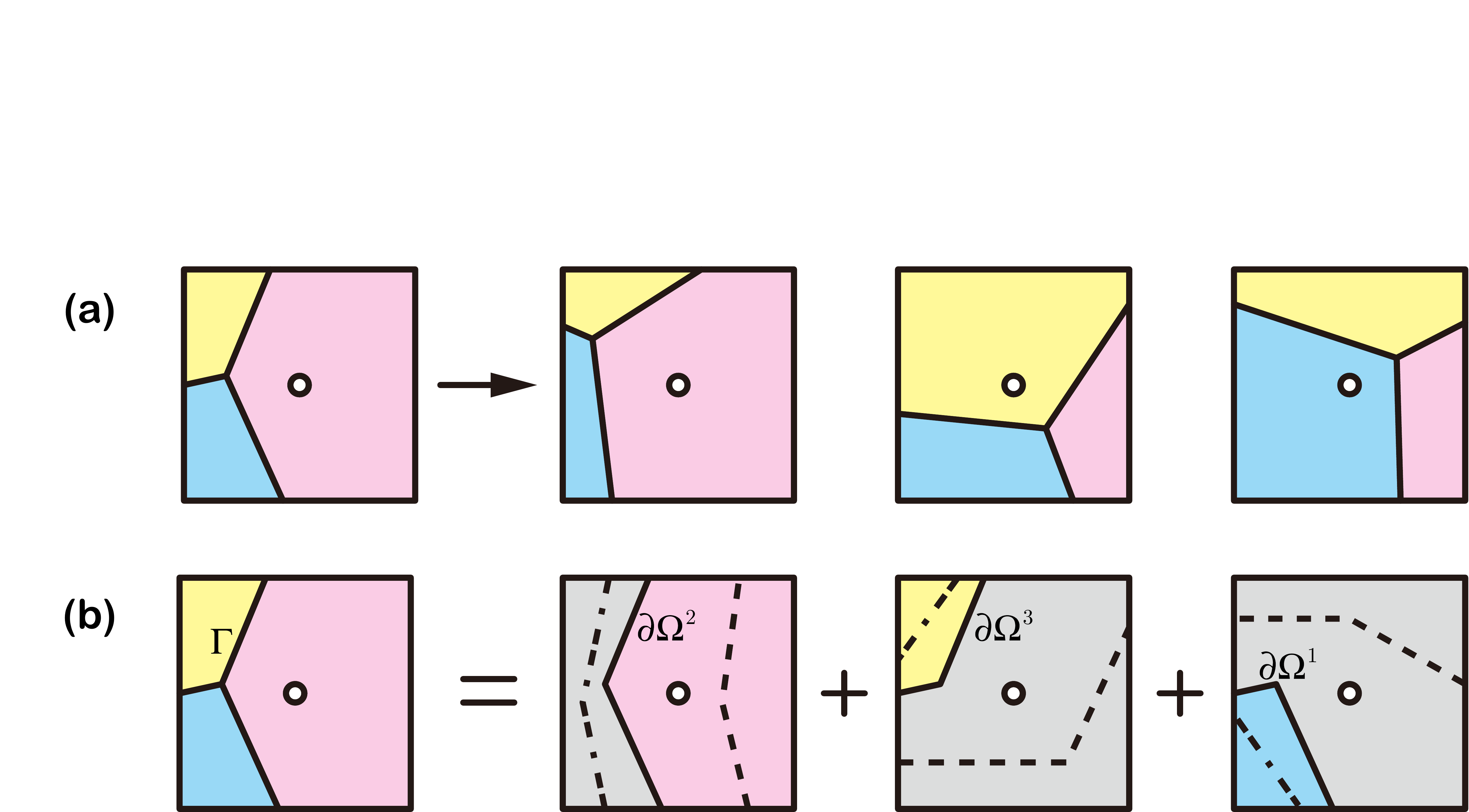}
\caption{A schematic representations of a 3-region (region domains are colored by blue($\chi=1$), yellow($\chi=3$) and red($\chi=2$)) cell.
(a) The evolution of this complex-region cut cell consisting of 3 regions, $\mathcal{N}_s = 3$.
Initially the primary indicator is $\chi_1=\chi_{i,j} = 2$,
and the secondary indicators are $\chi_2 = 1$ and $\chi_3 = 3$.
After one sub-step of advection, the primary indicator may be unchanged,
or have changed to $\chi_1 = 1$ or $\chi_1 = 3$.
(b) the surface evolution for each individual region.}
\label{stencil}
\end{center}
\end{figure}

\begin{figure}[p]
\begin{center}
\includegraphics[width=1.0\textwidth]{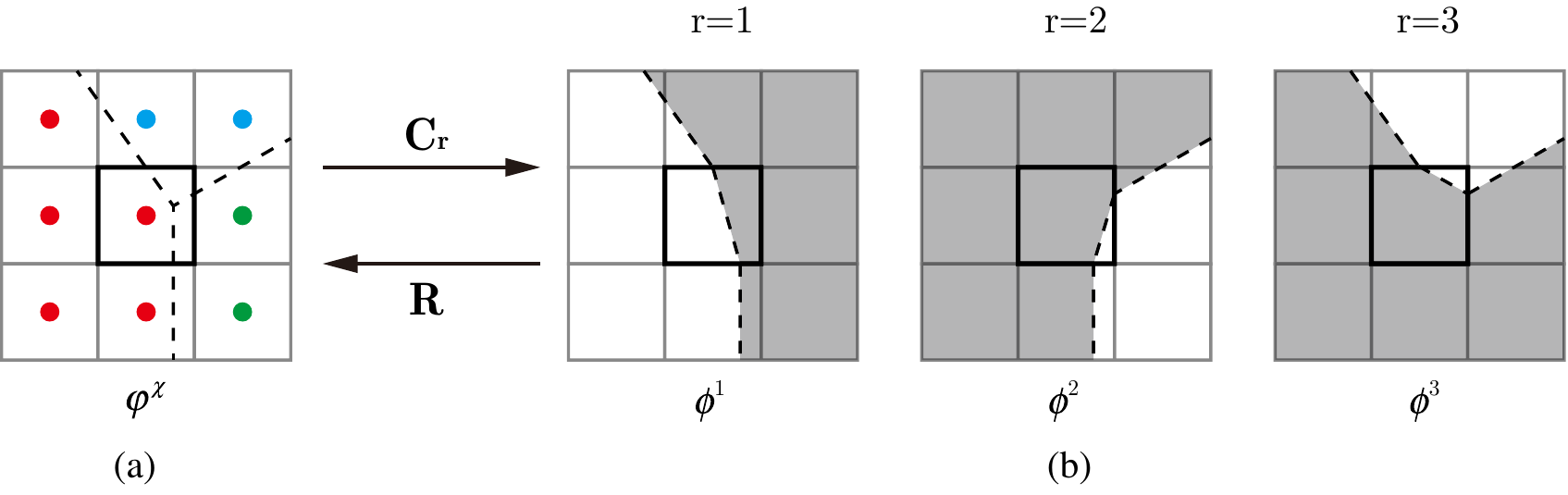}
\caption{A schematic representation example of local construction and reconstruction operators. Globally the system has $5$ regions and the local index set in this $3\times3$ stencil is $\mathrm{X}_s=\{r|r=1,2,3\}$. The corresponding region index $\chi_r$ is $\{\chi_r| r \in \mathrm{X}_s\} = \{4,2,5\}$. The mapping from (a) to each field in (b) is defined as the construction operator $\mathbf{C}: \mathbb{R}\times \mathbb{N} \rightarrow \mathbb{R}$. The inverse mapping from all fields in (b) to (a) is defined as the reconstruction operator $\mathbf{R} : \mathbb{R}^{3} \rightarrow \mathbb{R}\times \mathbb{N}$.
(a) The regional level-set field is $\varphi^{\chi}$. The center of cells whose region indicator equals $4$, $2$, and $5$ is colored by red, green, and blue, respectively.
(b) The multiple local level-set functions are $\phi^1$, $\phi^2$ and $\phi^3$ after applying $\mathbf{C}$ on $\varphi^{\chi}$. The gray part and white part correspond to the negative and positive $\phi^r$, respectively. The dashed line in each local level-set field is the region boundary and is represented by the zero contour of $\phi^r$.
}
\label{operator}
\end{center}
\end{figure}

The evolution of the interface network in a multi-region system
is equivalent to that of a signed level-set field whose interface advection is determined by the advection equation
\begin{equation} \label{advection-equation}
\frac{\partial \phi }{\partial t} + \mathbf{v}\cdot \nabla \phi = 0,
\end{equation}
where $\mathbf{v}$ is the advection velocity. The advection equation in our method is formulated locally for each region, identified by the region indicator, and recovers the original level-set method for cells that are sufficiently far away from a multiple junction. More detailed considerations are necessary to predict the evolution of the interface for cells in which more than two regions meet. Consider the situation in Fig.\ref{stencil}. The cell contains three regions colored by red, blue, and yellow. Initially, the largest fraction of the cell is occupied by region $\Omega^{\chi_1}$, with $\chi_1$ being the primary indicator. After advection by one time step, the interface has three possible configurations, as illustrated in Fig.\ref{stencil}(a). The first is that the cell center still resides within $\Omega^{\chi_1}$, so that the primary indicator of the cell does not change. In the other two cases, the primary indicator changes to the secondary indicator, either $\chi_2$ or $\chi_3$. After solving an advection equation for the local level-set field of the three regions separately, each region boundary may shrink (dashed line) or expand (dash-dotted line), see Fig.\ref{stencil}(b). For example, when the boundary $\partial\Omega^{\chi_1}$ of $\Omega^{\chi_1}$ moves away from the cell center its local level-set $\phi^{\chi_1}$ value increases (positive). For the converse case, $\partial\Omega^{\chi_1}$ moves across the cell center, so that it does not belong to $\Omega^{\chi_1}$ any longer, resulting in a corresponding sign change to $\phi^{\chi_1} < 0$. That the cell center is not in $\Omega^{\chi_1}$ does not necessarily imply that it is in $\Omega^{\chi_2}$ or $\Omega^{\chi_3}$. A direct combination of all three independently advected region boundaries may introduce an overlap or a void. To address this issue, we apply the reconstruction operator $\mathbf{R}$ to determine the regional level set $\varphi^{\chi}$ from the three candidate local level-set fields $\{\phi^{\chi_1}, \phi^{\chi_2}, \phi^{\chi_3}\}$. This is physically reasonable because $\mathbf{R}$ identifies the most likely indicator corresponding to the region domain in which the current cell center will be located after one advection time step.

The update of the regional level set of each finite-volume cell
$\varphi^{\chi,n}_{i,j} = (\varphi_{i,j}, \chi_1)^{n}$ at time-step $n$,
by a sub-step of an explicit time-integration scheme,
such as a strongly stable Runge-Kutta scheme \cite{shu1988efficient},
consists of three sub-steps:
(1) construction of the local multiple signed level-set fields in the stencil neighborhood $V_{c}$;
(2) computation of the new intermediate multiple local level-set fields by updating the locally constructed advection equations, see below;
(3) reconstruction of the regional level set from the multiple local level-set fields.
The local construction of the advection equations depends on the type of $V_{s}$ (full, two-region or complex-region).
Note that we assume that under the Courant-Friedrichs-Lewy (CFL) condition, i.e. $\mathrm{CFL} \leq 1.0$,
the type of $V_{s}$ does not change and the region indicators within the square-shaped near neighborhood remain unchanged during
one sub-step of the time-integration scheme, whereas the actual finite-volume cell may change its type. Unlike the ``repairing'' (or ``modification'') procedure in Ref. \cite{merriman1994motion}, no \textit{a posteriori} operations are required after advection.

\subsubsection{Updating a cell with full type $V_{s}$}
For a cell with full type $V_{s}$,
the local level set $\phi^{1, (n)}_{i,j} = \phi^n_{i,j}(\chi_1)$
and the intermediate data $\phi^{1, (s)}_{k,l} = \phi^{(s)}_{k,l}(\chi_1)$ in the stencil neighborhood
are constructed referring to the primary indicator $\chi_1$.
The intermediate-step update $\phi^{1, (s+1)}_{i,j}$ is obtained from
\begin{eqnarray} \label{local-advection}
&\phi^{1, (s+1)}_{i,j} = \alpha_s\phi^{1, (s)}_{i,j}
+ (1-\alpha_s)\left[\phi^{1, (s)}_{i,j} -
\Delta t \mathbf{v}^{n}_{i,j}\cdot\left(\nabla \phi^{1}\right)^{(s)}_{i,j}\right],
\quad s = 0,\ldots, m, \nonumber \\
&\phi^{1, (0)}_{i,j} = \phi^{1, n}_{i,j}, \quad \phi^{1, (m+1)}_{i,j} = \phi^{1, n+1}_{i,j}
\end{eqnarray}
where $m$ is the number of sub-steps, $\alpha_s$ is the parameter of Runge-Kutta sub-step $s$,
and $\left(\nabla \phi^{1}\right)^{(s)}_{i,j}$
is the finite difference approximation of the spatial derivative at the center of a finite-volume cell $C_{i,j}$.
The updated regional level set at the sub-step $s+1$ is
$\varphi^{\chi,(s+1)}_{i,j} = (\phi^{1, (s+1)}_{i,j}, \chi_1)$.

\subsubsection{Updating a cell with two-region type $V_{s}$}
For a cell with two-region type $V_{s}$,
the construction of the multiple local level-set fields and the regional level-set field updates are obtained essentially by the same operations as for a full cell.
The difference is that the intermediate regional level-set field update depends on the sign of $\phi^{1, (s+1)}_{i,j}$
\begin{equation} \label{two-region-cell}
\varphi^{\chi, (s+1)}_{i,j} = \begin{cases}
(\phi^{1, (s+1)}_{i,j}, \chi_1) \quad \rm{if}~ \phi^{1, (s+1)}_{i,j} \ge 0 \cr
(-\phi^{1, (s+1)}_{i,j}, \chi_2) \quad \rm{otherwise}	
\end{cases}.
\end{equation}
\subsubsection{Updating a cell with complex-region type $V_{s}$}
For a cell with complex-region type $V_{s}$,
the local level-set field is constructed and updated referring to the primary indicator $\chi_1$.
If $\phi^{1, (s+1)}_{i,j}$ does not change sign,
the new intermediate local level set is $\varphi^{\chi,(s+1)}_{i,j} = (\phi^{1, (s+1)}_{i,j}, \chi_1)$.
Otherwise, the local level-set field is constructed and updated
referring to the secondary indicators 
\begin{equation} \label{complex-region-cell}
\phi^{m, (s+1)}_{i,j} = \alpha_s\phi^{m, (s)}_{i,j}
+ (1-\alpha_s)\left[\phi^{m, (s)}_{i,j} -
\Delta t \mathbf{v}^{n}_{i,j}\cdot\left(\nabla \phi^{m}\right)^{(s)}_{i,j}\right],
\quad  \chi_m \in  V_{s},
\end{equation}
and the regional level set is reconstructed by the operator $\mathbf{R}$.

The reason for using different advection strategies based on the cell types is the following. A direct transport of the global unsigned level-set function which exhibits a discontinuity across the interface requires numerical diffusion for stabilization and thus produces a smeared interface. Consequently, the unsigned scalar function does not maintain the distance function property. Advection with a high-order, low dissipation scheme of the local signed level-set functions, however, maintains the sharp interface and the distance property much more accurately. We advect the unsigned level-set function for full type cells and the local constructed signed level-set functions for other cells. This essentially is an adaptive algorithm which first determines the cell type as the indicator of the local smoothness of the unsigned level-set. Then the unsigned level-set field is advected wherever the unsigned level-set is smooth or the signed level-set field is advected wherever the unsigned level-set is singular. The re-initialization is not subject to regularization constraints, and our numerical examples show that it can be chosen as infrequent as with the original level-set method.

Computational efficiency can be further increased by employing the narrow band technique \cite{adalsteinsson1995}. For simple test cases, such as in Sec. \ref{sec:triple_advect} and Sec. \ref{sec:scheme}, this is not necessary. For more complex cases, such as in Sec. \ref{sec:normal_flow} and Sec. \ref{sec:foam}, it is applied.
The operation count of our method per time step is $\mathcal{O}(n^2)$ (same as with the Semi-Lagrangian regional level-set method and original level-set method) and can be reduced to $\mathcal{O}(kn)$ by employing the narrow band technique \cite{adalsteinsson1995}, where $k$ is the band width and $n$ is the number of cells in any direction. For the entire system the construction and reconstruction operators are applied $d(2t+1)N_1 + 2d(2t+1)N_2 + \mathcal{N}_s d(2t+1)N_3$ and $N_3$ times respectively, where $N_1$, $N_2$ and $N_3$ are the number of cells with full type, two-region type and complex-region type $V_{s}$ in the narrow band of interface network, respectively. At every Runge-Kutta sub-step we solve the local advection equations $N_1+N_2+\mathcal{N}_s N_3$ times.

\subsection{Re-initialization}
It is important to re-initialize the regional level set when necessary for maintaining its distance-function property with respect to the interface network. We employ two re-initialization methods. One widely used method is that of Sussman et al. \cite{sussman1994level} where a re-initialization equation is solved iteratively until steady-state is reached. The implementation for the regional level-set method involves the following steps region by region. For each subdomain $\Omega^{\chi}$, $\chi \in \mathrm{X}$, first the local level set $\phi^r$ is constructed and the re-initialization equation according to Ref. \cite{sussman1994level} is iterated until a steady state is reached.
Alternatively, the explicit one-step method developed by Fu et al.\cite{fu2015} can be used which is significantly faster than iterative methods.

\subsection{Summary of the numerical method} \label{sec:summary}
Here we summarize our high-resolution regional level-set method and comment on the description above. The entire numerical method contains the following steps:

1. \textbf{The initialization step.} According to the description in Sec. \ref{sec3:0}, we define an initial regional level set $\varphi^{\chi}_{i,j} = (\varphi_{i,j}, {\chi}_{i,j})$ at each finite-volume cell $C_{i,j}$ center $\mathbf{x}_{i,j}$, where $\varphi_{i,j}$ is the distance from the interface network and ${\chi}_{i,j}$ indicates which region domain the center of $C_{i,j}$ is located at.\\
2. \textbf{The evolution step} contains three sub-steps for $C_{i,j}$: (1) construction of the $\mathcal{N}_s$ local level-set fields $\phi^{r,n}_{k,l}$ for the current time-step $n$ at the center of each cell $C_{k,l}$ which belongs to the stencil $V_c$ of $C_{i,j}$; (2) computation of $\phi^{r,n+1}_{i,j}$ at the next time-step $n+1$ by solving the $\mathcal{N}_s$ local advection equations; (3) reconstruction of the new regional level set $\varphi^{\chi}_{i,j}$ at the center of $C_{i,j}$ from the $\mathcal{N}_s$ new local level-set fields $\phi^{r,n+1}_{i,j}$ by the reconstruction operator $\mathbf{R}$. \\
3. \textbf{The re-initialization step.} Enforce the distance function property of the regional level set $\varphi^{\chi}_{i,j}$ as the distance from the interface network.We emphasize that the interface-network transport re-initialization not always is necessary, see Sec. \ref{sec4}.\\
4. \textbf{The postprocessing step.} If necessary, extract the interface network by the triangulation method of Ref. \cite{saye2012analysis}.  \\

It should be mentioned that except for a postprocessing step the method does not need to extract explicitly the interface network, so it is a fully implicit method for multi-region problems. Another important feature of our method is that unlike using neighbouring $\epsilon$-level-set contours to reconstruct the interface network \cite{saye2011voronoi}, our method directly captures the evolution of the interface network. The method is not a hybrid of the multiple level-set method and the regional level-set method. It rather can be viewed as a regional level-set method employing locally signed level-set fields.

\section{Numerical validation} \label{sec4}
In this section, we assess accuracy and efficiency, the present method by a range of numerical examples. Our intention is to show that we recover the high computational efficiency of the Semi-Lagrangian (SL) method, improves, however, significantly the prediction accuracy.
First, we compare different high-order discretizations of the level-set transport equation with SL results for constant rotation motion of a three-region case. Afterwards, two simple mean curvature flows are considered to verify the construction operators and re-initialization methods. Subsequently, the single vortex flow is used to demonstrate the ability of the present method to resolve long thin filaments. The computations of normal driven flow, mean curvature flow and their combination serve to assess the accuracy of the present method. Finally, we couple the present method with Navier-Stokes equations applied to dry-foam dynamics which undergo sudden breakups, in order to demonstrate the capability for coping with complex configurations.
\subsection{Simple test cases}
Three simple cases are considered to test suitable high-resolution finite-difference schemes, construction operators, and re-initialization methods for the simulation of multi-region problems. We consider the 5th-order weighted essentially non-oscillatory scheme (WENO) \cite{jiang1995efficient} and the central-upwind weighted essentially non-oscillatory scheme (WENO$\_$CU6) \cite{hu2010adaptive}. Both are compared with the SL scheme of the regional level-set method \cite{zheng2009simulation}.
\subsubsection{Circle expansion} \label{sec:expand}
\begin{figure}[p]
\begin{center}
\includegraphics[width=1.0\textwidth]{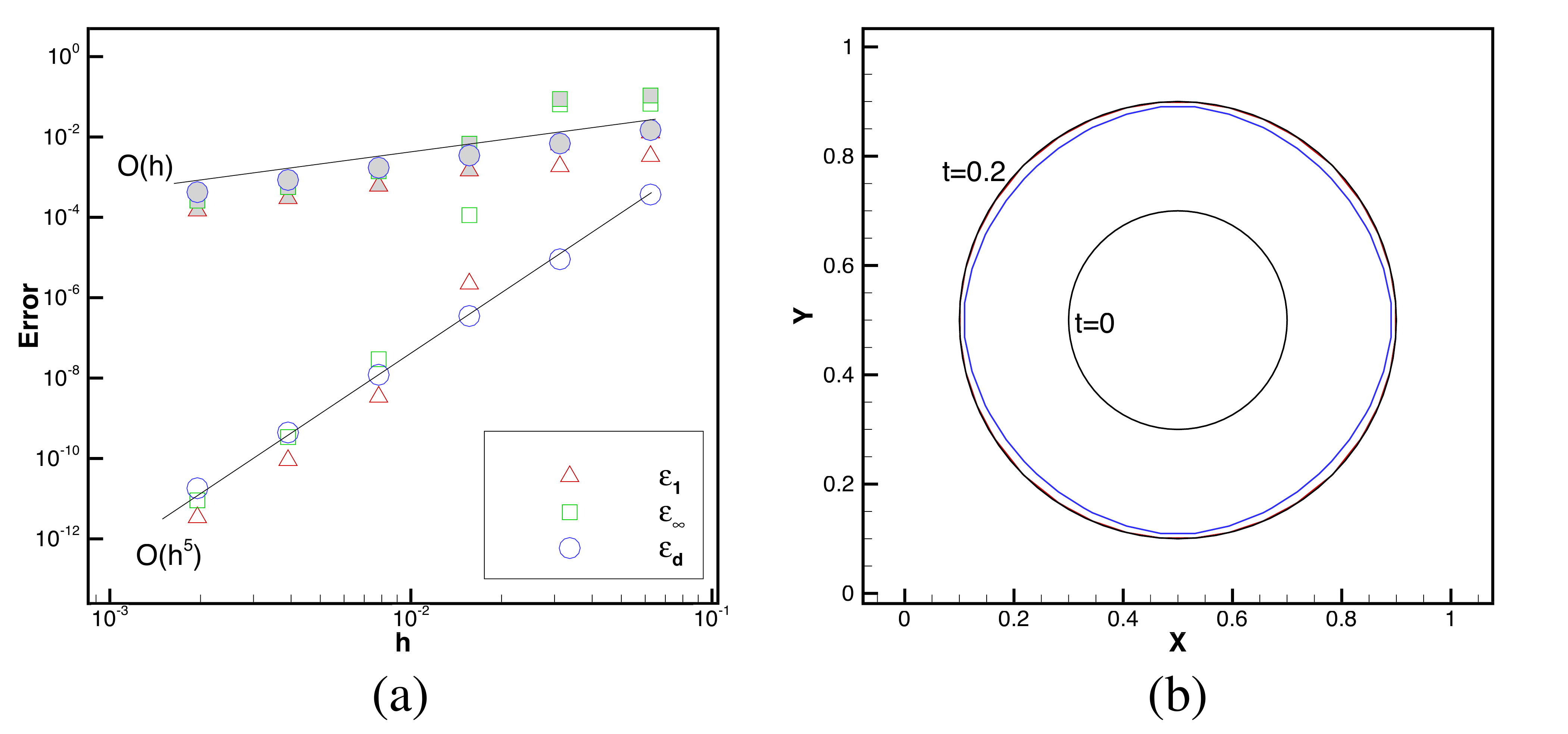}
\caption{Errors $\varepsilon_1$, $\varepsilon_{\infty}$ and $\varepsilon_d$ with increasing resolution (a) and the interfaces segmentation ($16 \times 16$ grid points) at $t=0.2$ (b) for a circle expansion. Present method (empty symbols and red line) is compared with Semi-Lagrangian regional level-set method (gray symbols and blue line). The exact interfaces at $t=0$ and $t=0.2$ are plotted with solid lines.}
\label{err_ana0}
\end{center}
\end{figure}

We start with a 2-region expansion case in Ref. \cite{saye2012analysis}. A circle with radius of $0.2$ expands with a uniform speed to $t=0.2$. The computational domain is $\left[0, 1\right]\times\left[0, 1\right]$. The explicit Euler scheme is used for time marching with a CFL number of $0.5$.The 5th-order WENO scheme is employed for spatial discretization. Fig. \ref{err_ana0}(a) shows three error measures
\begin{equation} \label{error}
\varepsilon_1 = \int_0^T \lVert \varphi \rVert_1\;\mathrm{d}x, \quad \varepsilon_{\infty} = \int_0^T \lVert \varphi \rVert_{\infty}\;\mathrm{d}x, \quad \varepsilon_d = \int_0^T d_H (\Gamma^n, \Gamma^e)\;\mathrm{d}x,
\end{equation}
with
\begin{eqnarray}
\lVert \varphi \rVert_1 &=&
\max_{\substack{(i,j) \in B}}| \varphi^n_{i,j}- \varphi^e_{i,j} |,  \qquad
\lVert \varphi \rVert_{\infty} =
\frac{1}{n(B)}\sum_{\substack{(i,j) \in B}}| \varphi^n_{i,j}- \varphi^e_{i,j} |, \nonumber \\
B&=&\{(i,j) | \varphi^n_{i,j} < 10\Delta x  \}.
\end{eqnarray}
The Hausdorff distance is
\begin{equation}
d_H (\Gamma^n, \Gamma^e)= \max( \sup\limits_{\textbf{y} \in \Gamma^e} \inf \limits_{\textbf{x} \in \Gamma^n} \lVert \textbf{x} - \textbf{y} \rVert_2, \sup\limits_{\textbf{x} \in \Gamma^e} \inf \limits_{\textbf{y} \in \Gamma^n} \lVert \textbf{x} - \textbf{y} \rVert_2).
\end{equation}
The superscripts `e' and `n' stand for the exact and the numerical solutions. Our method achieves the expected 5th-order convergence rate for all three error measures. The SL method is 1st-order, as expected, as is the Voronoi implicit interface method, see Fig. 4 of Ref. \cite{saye2012analysis}. Fig. \ref{err_ana0}(b) shows the interfaces of our method and SL method at $t=0.2$, compared the exact solution.
\subsubsection{Triple point advection} \label{sec:triple_advect}
\begin{figure}[p]
\begin{center}
\includegraphics[width=1.0\textwidth]{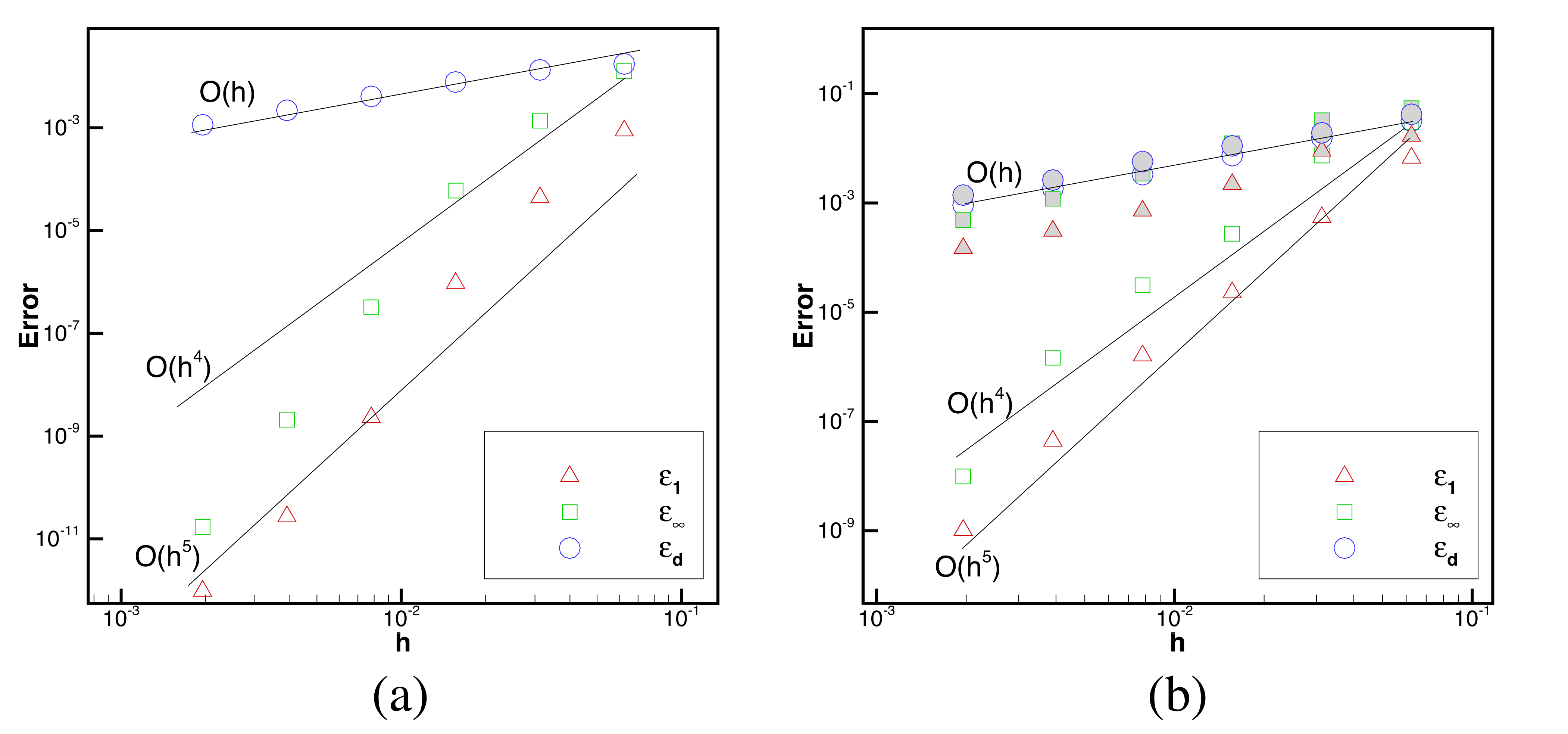}
\caption{Errors $\varepsilon_1$, $\varepsilon_{\infty}$ and $\varepsilon_d$ with increasing resolution for the single (a) and double (b) triple point advection. Present method (empty symbols) is compared with SL regional level-set method (gray symbols).}
\label{err_ana}
\end{center}
\end{figure}

We assess the prediction accuracy of our method by triple point advection cases which have analytical solutions. The first configuration contains only one triple point which initially is located at $(0.2, 0.5)$. The second is a circle of radius $r_0=0.3$ divided into two parts. The computational domain is $\left[0, 1\right]\times\left[0, 1\right]$. The velocity field is given by $(u,v) = (1.0, 0.0)$. The explicit Euler scheme is used for time marching with a CFL number of $0.6$. The 5th-order WENO scheme is used for spatial discretization. The simulations are performed until time $t=0.4$. Re-initialization is not employed. Error measures of Eq. (\ref{error}) are computed within a narrow band
\begin{eqnarray}
B=\{(i,j) | \varphi^n_{i,j} < 1.2\Delta x  \cap  \lVert \textbf{x}_{i,j} - \textbf{x}_s  \rVert_2 < 0.05L  \},
\end{eqnarray}
where $\textbf{x}_{i,j}$ and $\textbf{x}_s$ are the locations of cell center and triple points, respectively. Due to the smoothness properties of the exact solution at most first order accuracy in $\varepsilon_d$ can be achieved. This is reproduced by our method, as shown in Figs. \ref{err_ana}(a) and \ref{err_ana}(b). The benefit of the high-resolution discretization becomes evident for the global measures $\varepsilon_1$ and $\varepsilon_{\infty}$. The interface recovers the exact solution as the resolution increases.
A comparison with the SL regional level-set method, shown in Fig. \ref{err_ana}(b), demonstrates that the error magnitude of the SL regional level-set method is significantly larger than that of our method.

\subsubsection{Constant rotation of three regions} \label{sec:scheme}
\begin{figure}[p]
\begin{center}
\includegraphics[width=0.9\textwidth]{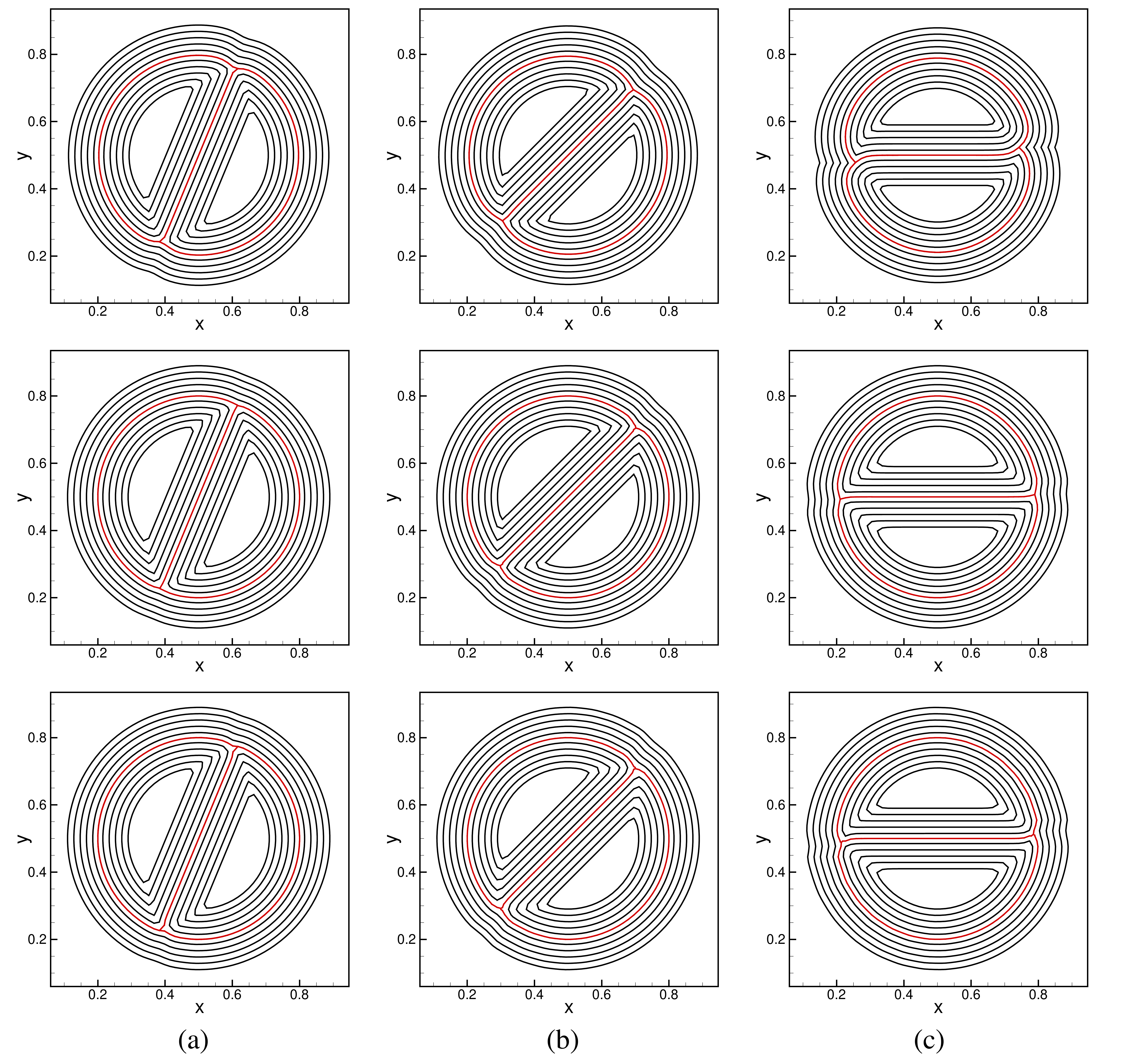}
\caption{Regional level-set contours ranging from 0.015 to 0.07: (a) $t=\frac{1}{8}\pi$, (b) $t=\frac{1}{4}\pi$, (c) $t=\frac{1}{2}\pi$. The results are obtained by using three different schemes: SL (first row), 5th-order WENO (middle row) and WENO$\_$CU6 (last row).}
\label{schemes}
\end{center}
\end{figure}

\begin{figure}[p]
\begin{center}
\includegraphics[width=1.0\textwidth]{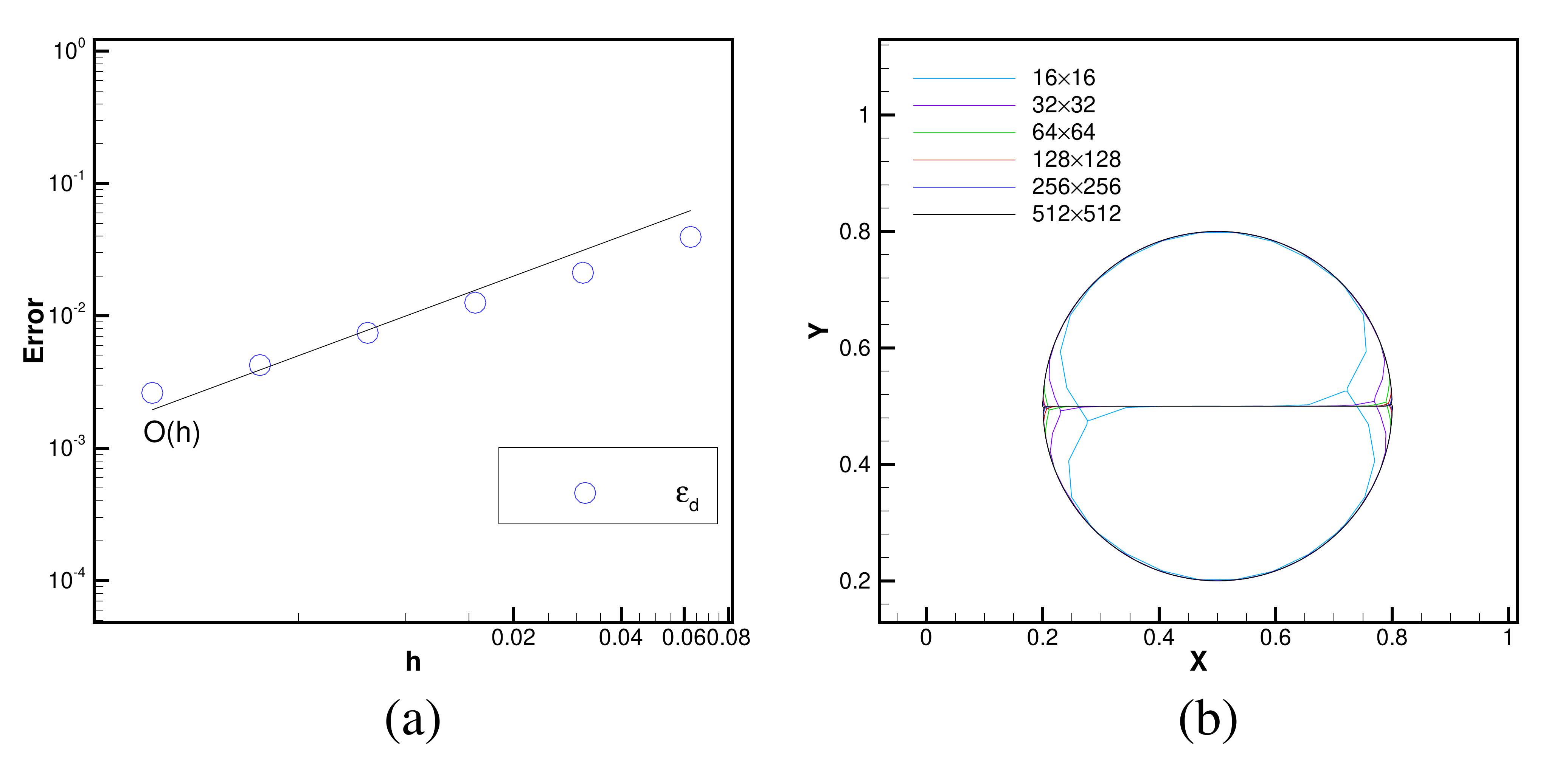}
\caption{(a) Error $\varepsilon_d$ with increasing grid resolution for constant rotation case. (b) Segmentation of the interface networks at $t=\frac{1}{2}\pi$.}
\label{err_ana2}
\end{center}
\end{figure}

A two-dimensional circle of radius $r_0=0.3$ is divided into two equal parts which undergo constant rotation. The computational domain is $\left[0, 1\right]\times\left[0, 1\right]$. Symmetry conditions are employed at all boundaries. The grid spacing is $h=\frac{1}{64}$. The explicit Euler scheme is used for time marching with a CFL number of 0.6. Two high-resolution schemes (5th-order WENO and WENO$\_$CU6), see the last two rows of Fig.\ref{schemes}, are compared with respect to their ability to capture the two triple points, in comparison with the results of the 1st-order SL scheme, see the first row of Fig.\ref{schemes}. Three regional level-set contours are shown in Fig.\ref{schemes} at $t=\frac{1}{8}\pi$, $t=\frac{1}{4}\pi$, and $t=\frac{1}{2}\pi$. Apparently, in each simulation, the deflection angles of the separation interface at all times agree with the theoretical values. However, the level-set contours surrounding the interface alter slightly. The SL method results in more smeared contours near the triple points due to numerical dissipation. The two high-resolution schemes have less numerical dissipation and thus preserve sharp corners. The relative radius difference $|\triangle r|/r_0$ of the circle for the WENO and WENO$\_$CU6 results at $t=0.5\pi$ are $0.211\%$ and $0.285\%$, respectively, and are much smaller than for the SL result which is $3.42\%$, indicating better area conservation with WENO and WENO$\_$CU6. An error analysis for this case, shown in Fig. \ref{err_ana2}(a), indicates 1st-order convergence of the interface location. The extracted interface network converges to the exact solution asymptotically as the resolution increases, as shown in Fig. \ref{err_ana2}(b).

\subsubsection{\label{sec:2dmean}Two dimensional mean curvature flows}
\begin{figure}[p]
\begin{center}
\includegraphics[width=1.0\textwidth]{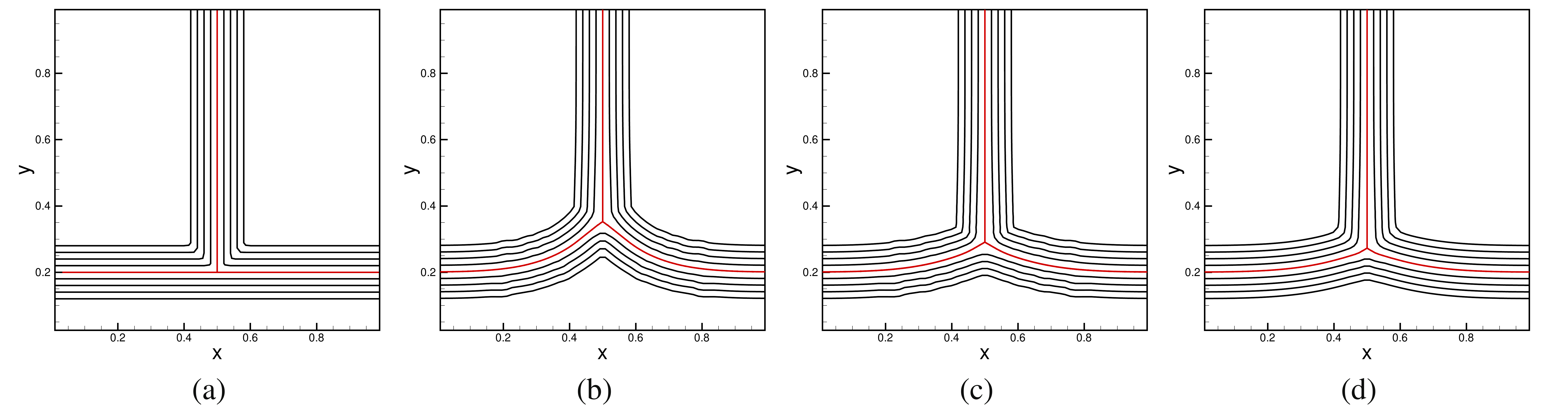}
\caption{Regional level-set contours ranging from 0.02 to 0.06: (a) initial contours; (b) t=0.02, explicit re-initialization method and $\mathbf{C}_r$; (c) t=0.02, explicit re-initialization method and $\mathbf{C}_r^{\star}$; (d) t=0.02, iterative re-initialization method and $\mathbf{C}_r^{\star}$.}
\label{Tjunction}
\end{center}
\end{figure}

The interface network of three-region and five-region systems is evolved under mean curvature flow with $\mathbf{u}=\kappa\mathbf{n}$, where $\mathbf{n}$ and $\kappa$ are calculated by Eqs. (\ref{normal}) and (\ref{curvature}). In our simulation we employ the 5th-order WENO scheme for advection and consider different construction operators and re-initialization methods. Since the mean curvature uses 2nd-order derivatives of $\phi^r$, it is more sensitive to the local constructed level-set fields. The computational domain extent is $\left[0, 1\right]\times\left[0, 1\right]$ in x and y directions. Symmetry boundary conditions are employed at all domain boundaries. The explicit Euler scheme is used for time marching, and the time step $\Delta t$ is the same as that in \cite{saye2012analysis}, $\Delta t = \frac{h^2}{4}$. As depicted in Fig.\ref{Tjunction}, the T junction transforms to a Y junction under the effect of mean curvature. The different columns in this figure show results for different re-initialization methods and construction operators. When we use the construction operator $\mathbf{C}_r$, the displacement of the triple point is overestimated in both cases, see Fig.\ref{Tjunction}(b). Significant improvement is observed by the operator $\mathbf{C}_r^{\star}$ where $\chi(\mathbf{x}) \neq \chi_r$, irrespective of the employed re-initialization method, see Figs. \ref{Tjunction}(c) and \ref{Tjunction}(d). Note that noise in the contours far away from the interface visible for the explicit re-initialization is irrelevant to our method as only the smooth inner contours are used for extracting the interface.

\begin{figure}[p]
\subfigure[] {\label{accuracy2da} \includegraphics[scale=0.3]{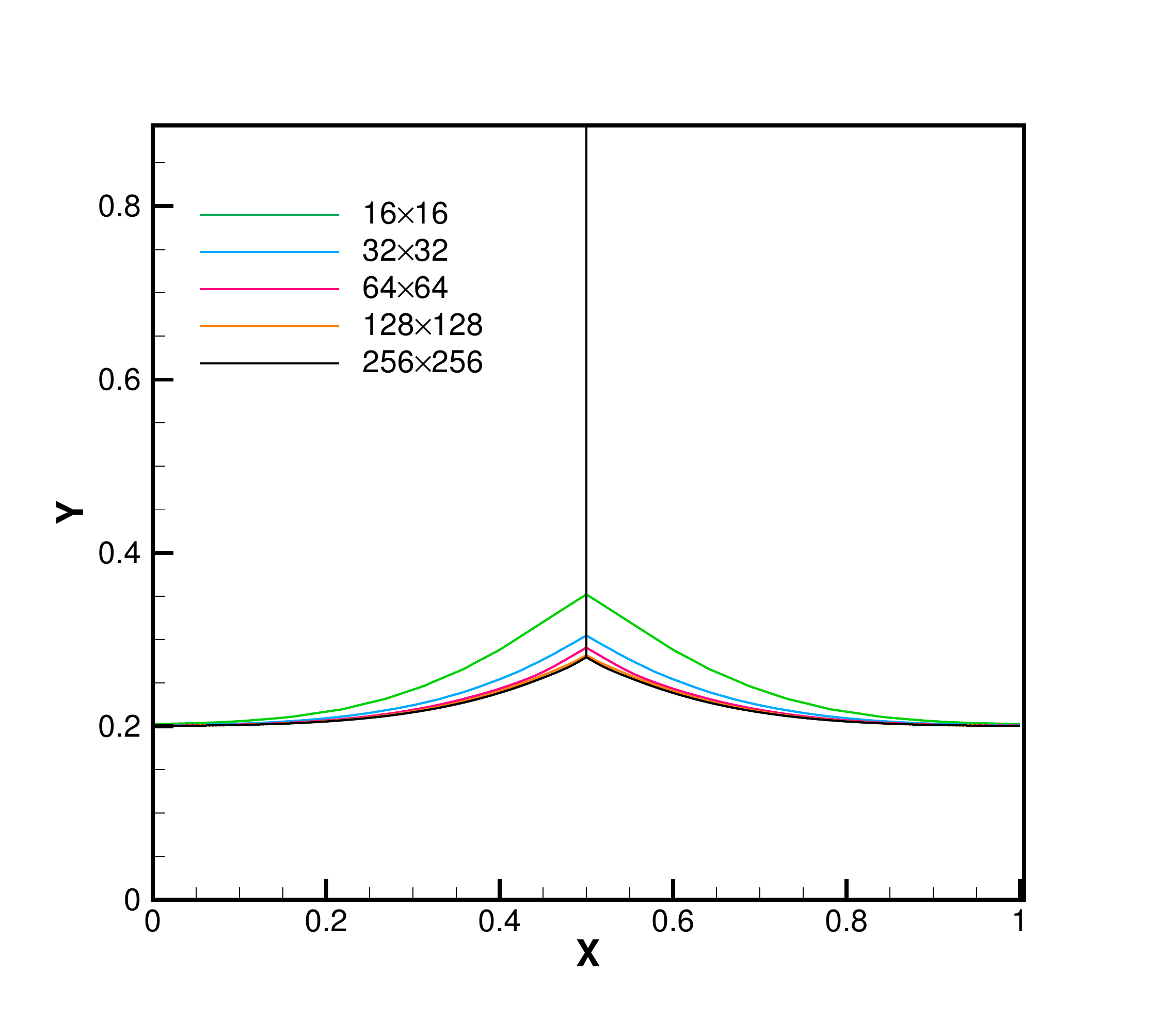}}
\subfigure[] {\label{accuracy2db} \includegraphics[scale=0.3]{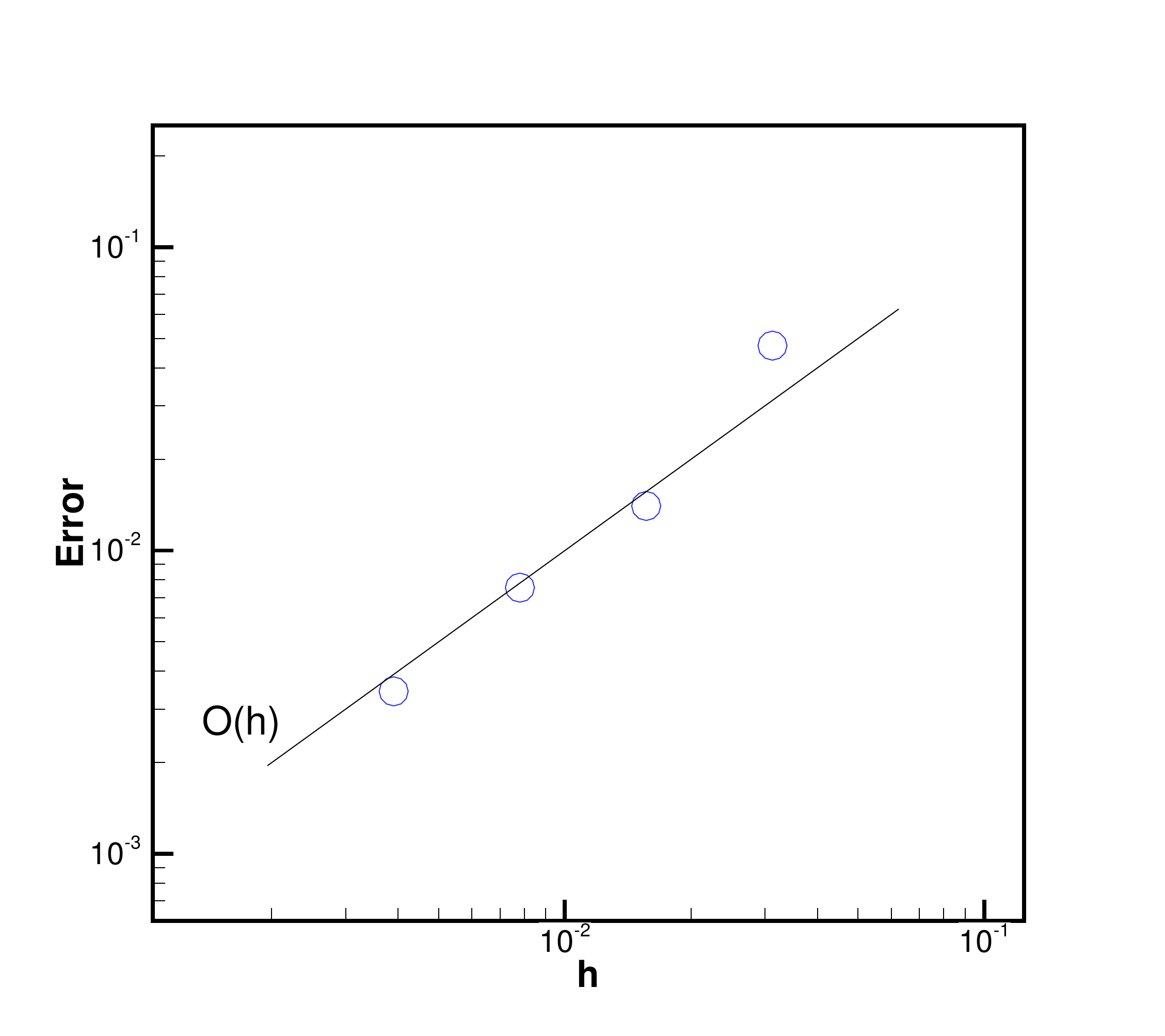}}
\caption{\label{accuracy2d} Convergence for mean curvature flows in Fig. \ref{Tjunction}: (a) interface segmentation and (b) error of triple point locations.}
\end{figure}

We conclude that the construction operator $\mathbf{C}_r^{\star}$ is more suitable for mean curvature flows. Both re-initialization methods are suitable for capturing the interface network. We obtain convergened interface locations, as shown in Fig. \ref{accuracy2d}(b), where the interface networks are extracted by the triangulation method of Ref. \cite{saye2012analysis}. As shown in Fig. \ref{accuracy2d}(b), first order convergence is achieved for the triple-point locations.

\subsection{Single vortex flow}
For the single vortex flow case, we use the setup introduced by Bell et al. \cite{bell1989second} to test the ability and accuracy on resolving thin filaments under the deformation by the velocity field:
\begin{eqnarray}
\left\{
\begin{aligned}
&u=-2{{\sin }^{2}}\left( \pi x \right)\sin \left( 2\pi y \right) \\
&v=2{{\sin }^{2}}\left( \pi y \right)\sin \left( 2\pi x \right)
\end{aligned}
\right. .
\end{eqnarray}
Again, we employ the 5th-order WENO scheme for advection. The initial circle deforms into a filament wrapping around the center of the domain. This structure wraps back into the initial circle upon reversing the velocity at $t=3$. The centers of concentric circles are $(0.5, 0.75)$ and their radii are $0.08$ and $0.22$, respectively. Symmetric boundary conditions and a second-order strongly stable Runge-Kutta scheme \cite{domingues2008adaptive} with a CFL number of 0.6 are employed. The grid size is refined from $h=\frac{1}{256}$ to $h=\frac{1}{1024}$. 

\begin{figure}[h]
\begin{center}
\includegraphics[width=1.0\textwidth]{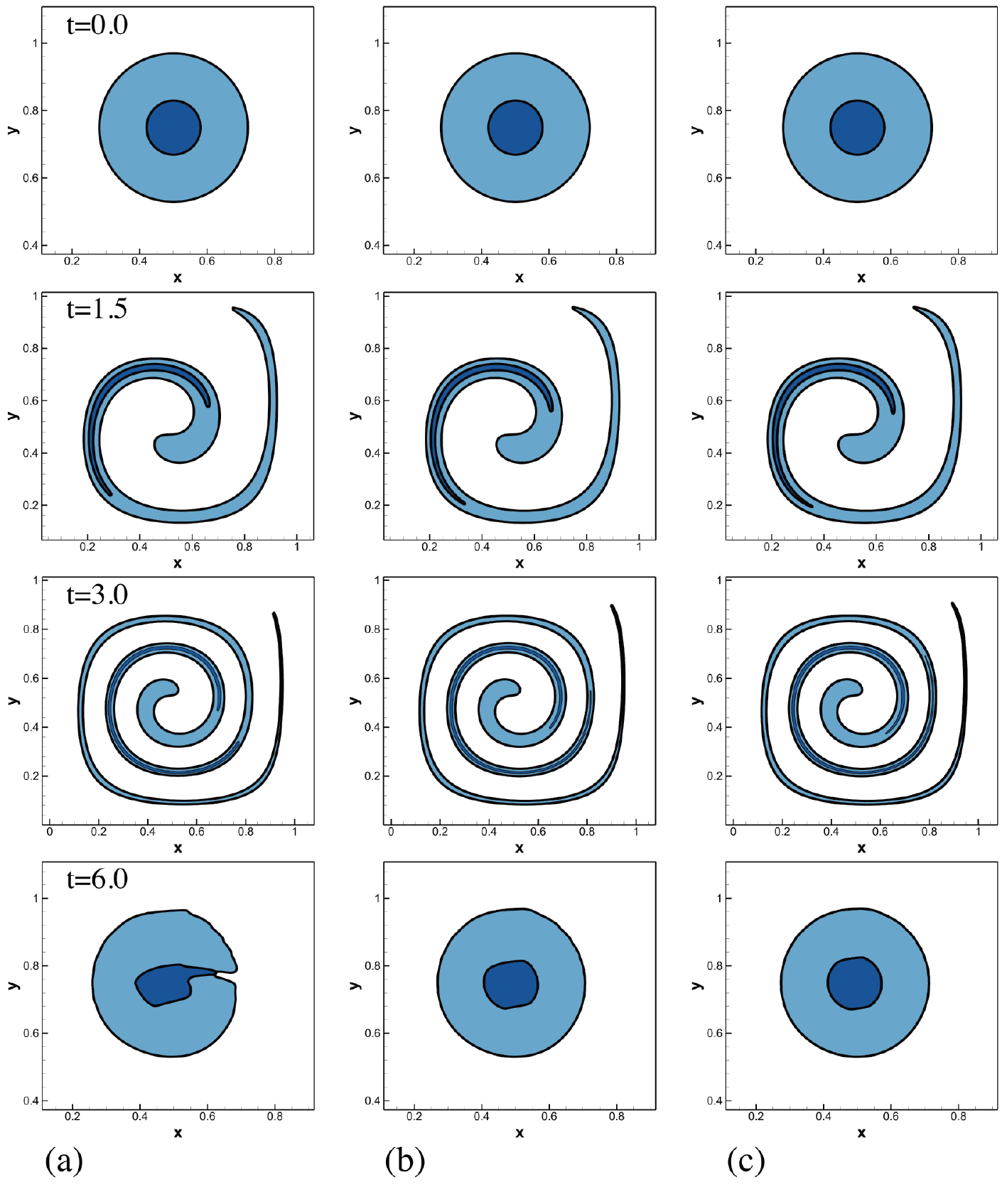}
\caption{Interface deformation for the single vortex flow on at $t = 0$, $1.5$, $3.0$, and $6.0$ with different resolutions: (a) $h=\frac{1}{256}$, (b) $h=\frac{1}{512}$, and (c) $h=\frac{1}{1024}$.}
\label{rotation}
\end{center}
\end{figure}

It can be seen from Fig.\ref{rotation} that two spirals of long filaments are successfully captured at $t=1.5$ and $3.0$. With increasing grid resolution the filaments become longer, especially the inner one. At $t=6.0$ with $h=\frac{1}{256}$ two additional unresolvable triple points are generated such that the two initial circle contours are connected. This phenomenon can be attributed to the fact that for extremely stretched filaments the two interfaces may reside in the same cell, and a small numerical perturbation leads to a topology change while reverse rotation does not disconnect the interface. We emphasize that although we use $256 \times 256$ cells overall, the number of cells across the initial inner circle is only $40 \times 40$. Increasing mesh resolution removes such artifacts, as shown in Figs.\ref{rotation}(b) and (c). Both circles recover their initial shapes with good area conservation. 

\begin{figure}[h]
\begin{center}
\includegraphics[width=1.0\textwidth]{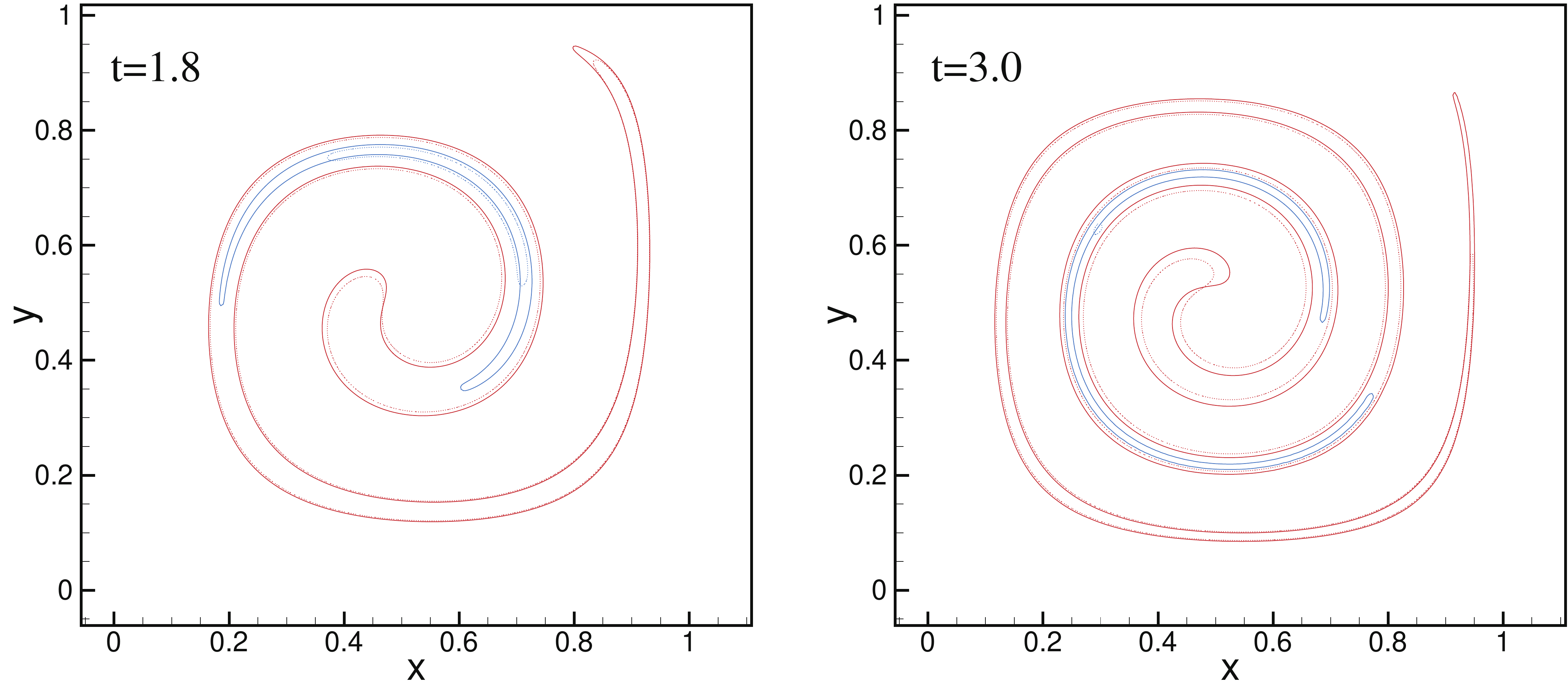}
\caption{Outer (red line) and inner (bule line) interfaces deformation for the single vortex flow with Semi-Lagrangian scheme (dotted line) and 5th-order WENO scheme (solid line).}
\label{rotation2}
\end{center}
\end{figure}

In Fig. \ref{rotation2}, we compare the results of our method with that of the SL regional level-set method. We observe that the SL method shows large numerical dissipation and can not reproduce the long filament. The area conservation errors are listed in Table. \ref{table1}. Note that the inner circle is very poorly resolved, i.e., for $h=\frac{1}{64}$, we have only $12 \times 12$ cells across the initial inner circle. With such low resolution $100\%$ mass loss of the inner circle is inevitable for any method.
At all three grid resolutions, our method exhibits significantly smaller area loss compared to the SL method. The CPU time measurement indicates that our method is nearly as fast as the SL regional level-set method even though we employ the high-order spatial schemes. 
We conclude that our approach achieves improved accuracy and comparable computational cost compared with the SL regional level-set method.
\begin{table}
\caption{\label{table1} Area conservation error ($\%$) and CPU times (in seconds) of the Semi-Lagrangian regional level-set method (SL-RLS) and our high-resolution regional level-set method (HR-RLS) for simulating the single vortex flow. The area loss $\Delta A(t_0) = (\Delta A^1(t_0), \Delta A^2(t_0))$ and $\Delta A(t_1) = (\Delta A^1(t_1), \Delta A^2(t_1))$ are measured at physical time $t_0=3.0$ and $t_1=6.0$, where the superscripts $1$ and $2$ indicate the outer and inner circles, respectively. The CPU time is measured at $t=3.0$.}
\resizebox{\columnwidth}{!}{
\begin{tabular}{lccccccccccc}
\hline\hline
\multirow{3}{*}{} Method& \multicolumn{3}{l}{$h=\frac{1}{64}$} & & \multicolumn{3}{l}{$h=\frac{1}{128}$} & &  \multicolumn{3}{l}{$h=\frac{1}{256}$} \\
\cline{2-4}\cline{6-8}\cline{10-12}
&$\Delta A(t_0)$ & $\Delta A(t_1)$ & CPU time &  & $\Delta A(t_0)$& $\Delta A(t_1)$& CPU time & &$\Delta A(t_0)$ &$\Delta A(t_1)$ &CPU time\\
\hline SL-RLS &(80.2, 100)	&(100, 100)	&5.32	&	&(24.9, 100)	&(51.9, 100)	&24.68	& &(7.6, 98.8) &(15.6, 100) &134.88  \\
\hline HR-RLS  &(72.4, 100)	&(72.8, 100) &6.44	&	&(24.5, 10.5)	&(41.0, 59.8)	&44.52	& &(2.3, 15.3) &(5.2, 17.6) &147.30 \\
\hline \hline
\end{tabular}
}
\end{table}
\vspace{4cm}

\subsection{Normal driven flow} \label{sec:normal_flow}
\begin{figure}[p]
\begin{center}
\includegraphics[width=1.0\textwidth]{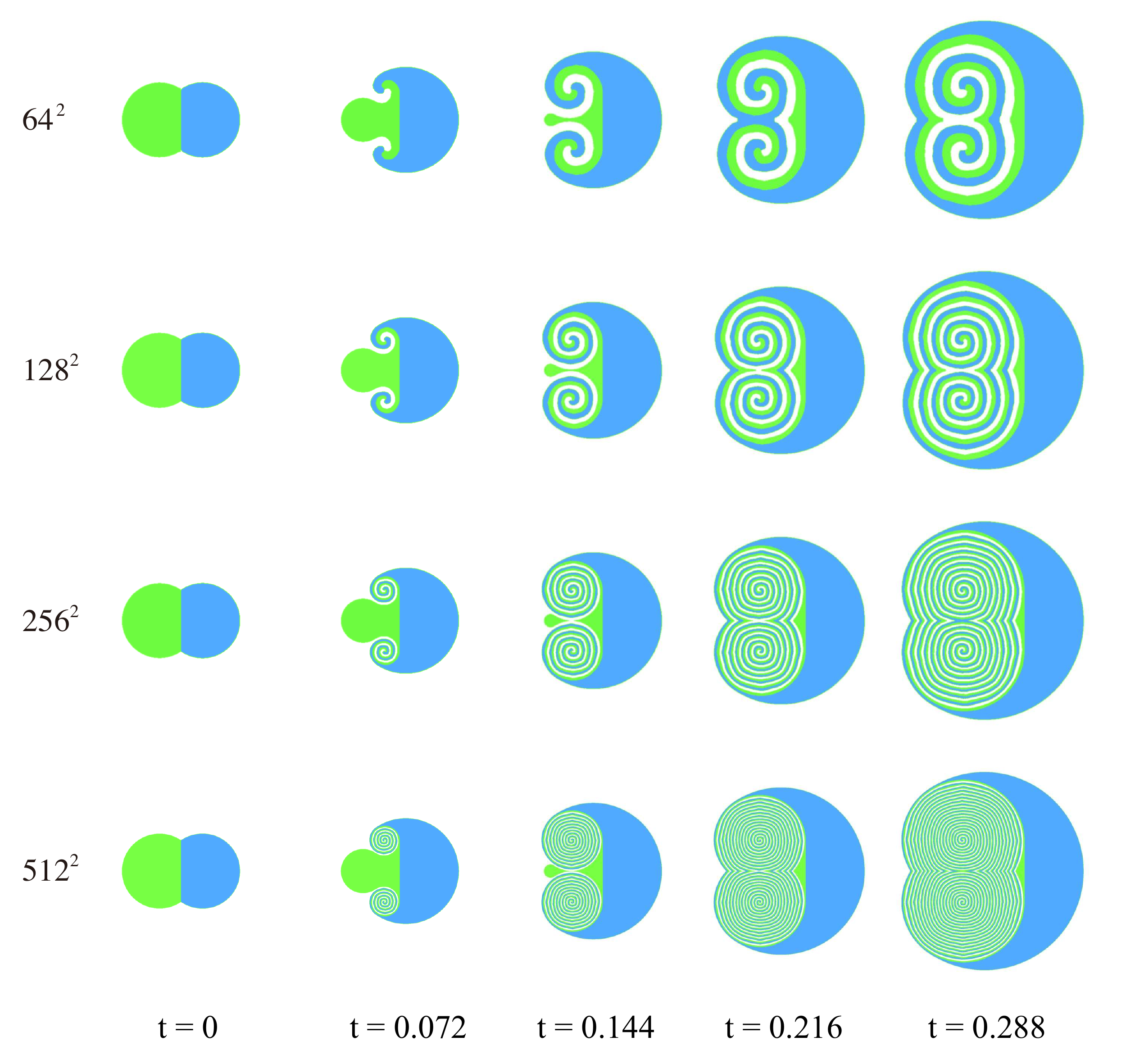}
\caption{Constant normal driven flow of three regions at $t = 0$, $0.072$, $0.144$, $0.216$, and $0.288$ with different resolutions ($h=\frac{1}{64}$, $h=\frac{1}{128}$, $h=\frac{1}{256}$, and $h=\frac{1}{512}$).}
\label{normalflow}
\end{center}
\end{figure}

The constant normal driven flow has been studied in \cite{saye2012analysis}. The pairwise $\Gamma_{ab}$ interface separating region domains $\Omega^a$ and $\Omega^b$ moves in its normal direction with a constant speed. As shown in Fig.\ref{normalflow}, two initially neighboring region domains $\Omega^b$ (colored by green) and $\Omega^c$ (colored by blue) are surrounded by a background region domain $\Omega^a$ (colored by white). We define the normal velocity of the interface of each region as
\begin{equation}
\mathbf{u}_{\Gamma_{ab}}= \mathbf{n}_{b}, \mathbf{u}_{\Gamma_{bc}}= \mathbf{n}_{c}, \mathbf{u}_{\Gamma_{ca}}= \mathbf{n}_{a},
\end{equation}
where $\mathbf{u}_{\Gamma_{ab}}$ is the velocity at the pairwise interface $\Gamma_{ab}$ and $\mathbf{n}_b$ is the normal direction of the region boundary $\partial \Omega_b$ at $\Gamma_{ab}$.
Thus the flow wraps $\Omega^{a}$ into $\Omega^{b}$ which in turn is wrapped into $\Omega^{c}$. The computation is carried out on a unit square with different resolutions employing the 5th-order WENO schemes for advection. Symmetry boundary conditions and the explicit Euler scheme are employed with a CFL number of $0.6$. 

\begin{figure}[p]
\begin{center}
\includegraphics[width=1.0\textwidth]{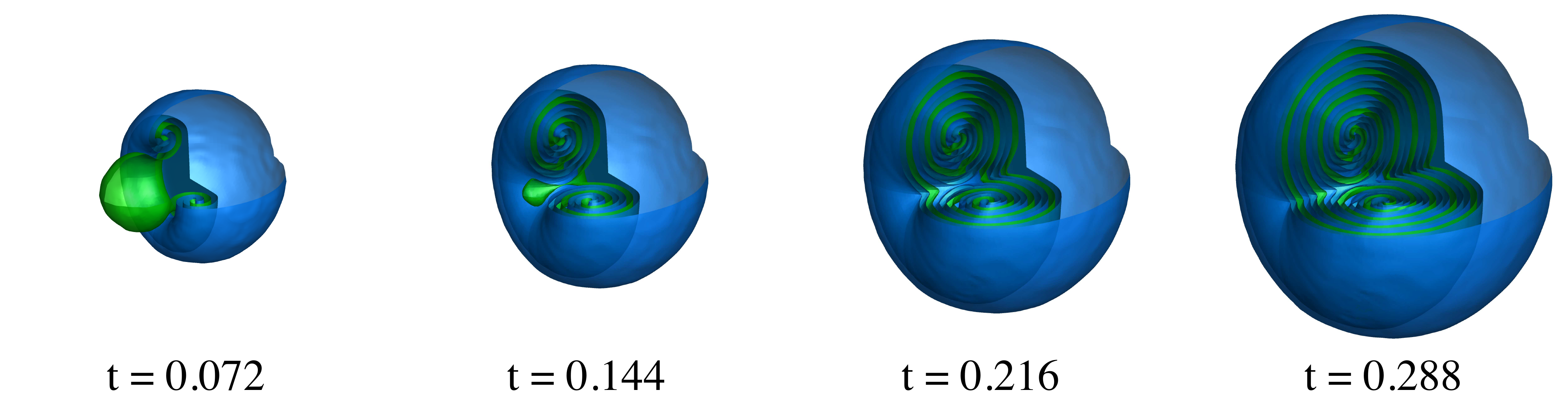}
\caption{3D normal driven flow at $t = 0.072$, $0.144$, $0.216$, and $0.288$ with $h=\frac{1}{128}$.}
\label{normalflow3D}
\end{center}
\end{figure}

\begin{figure}[p]
\subfigure[] {\label{spirala} \includegraphics[scale=0.3]{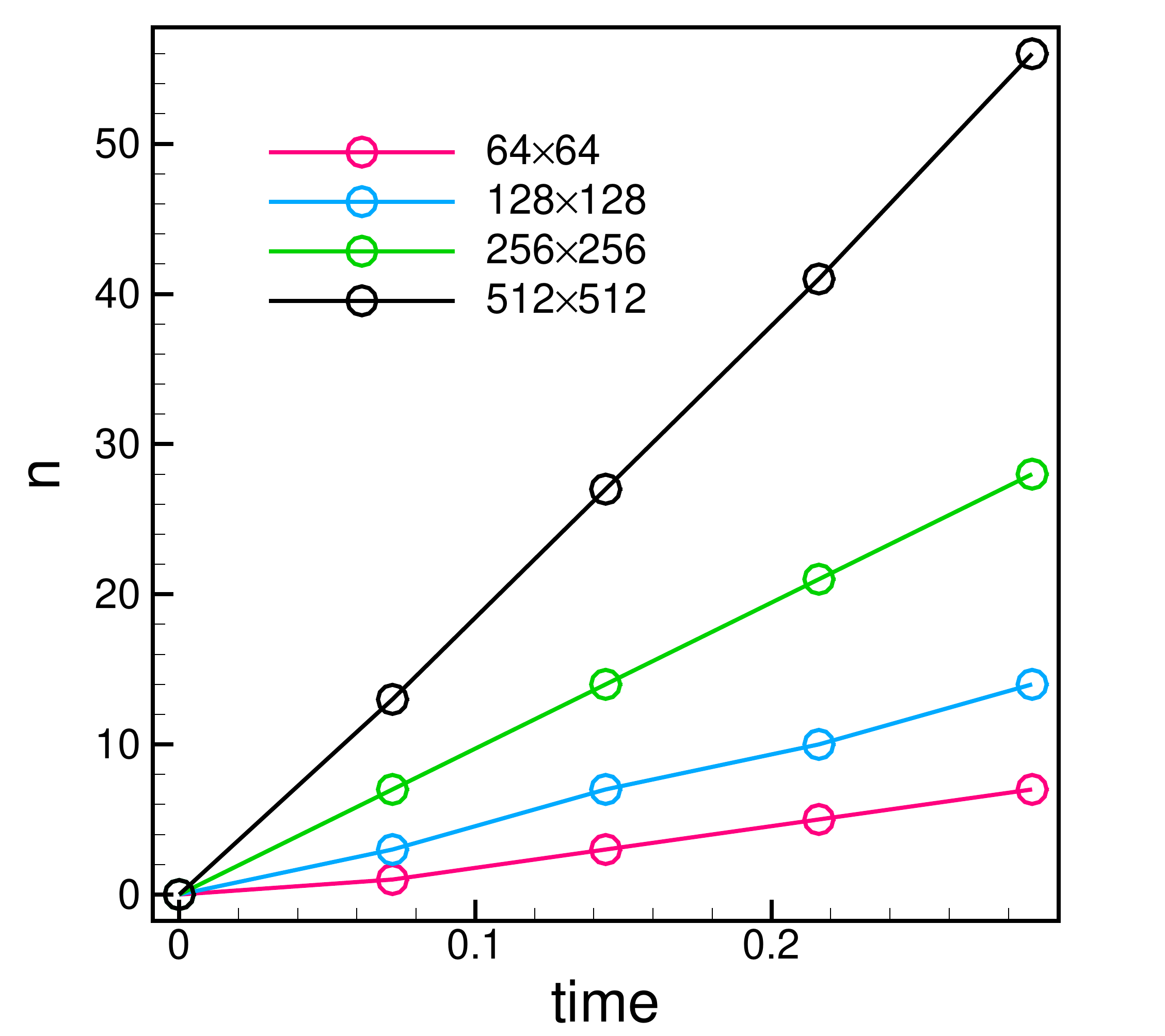}}
\subfigure[] {\label{spiralb} \includegraphics[scale=0.3]{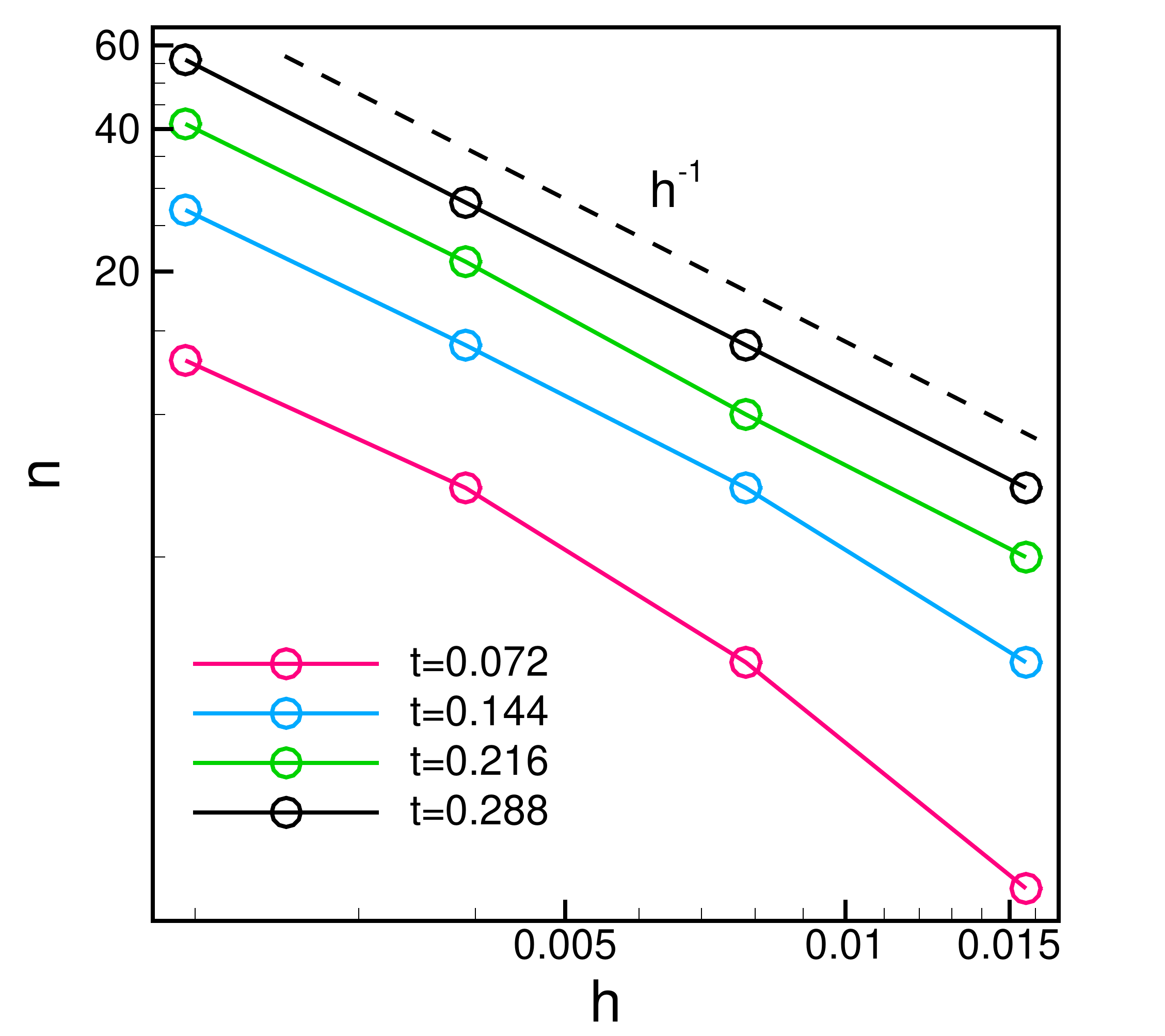}}
\caption{\label{spiral} Number of spirals for normal driven flow: (a) temporal evolution and (b) relation with resolution.}
\end{figure}

For a comparison between the results of our method and those in \cite{saye2012analysis} we plot five snapshots at the same time instants as those of \cite{saye2012analysis}, from $t=0$ to $t=0.288$. During evolution, the number of spirals increases quickly, and we observe that our results actually exhibit more visible spirals compared to that in \cite{saye2012analysis}, both on the same $256 \times 256$ grid, see third row of Fig.\ref{normalflow}. This observation can be attributed to the high-resolution scheme employed in our method. Another reason is that we directly advect the interface network.
The number of spirals increases linearly with time and proportionally to grid resolution, see Figs \ref{spirala} and \ref{spiralb}, respectively. In order to demonstrate the validity of our method in three dimensions, a 3D example of normal driven flow on a $128 \times 128 \times 128$ grid, is shown in Fig. \ref{normalflow3D}.

\subsection{Von Neumann-Mullins' law validation}
\begin{figure}[p]
\begin{center}
\includegraphics[width=1.0\textwidth]{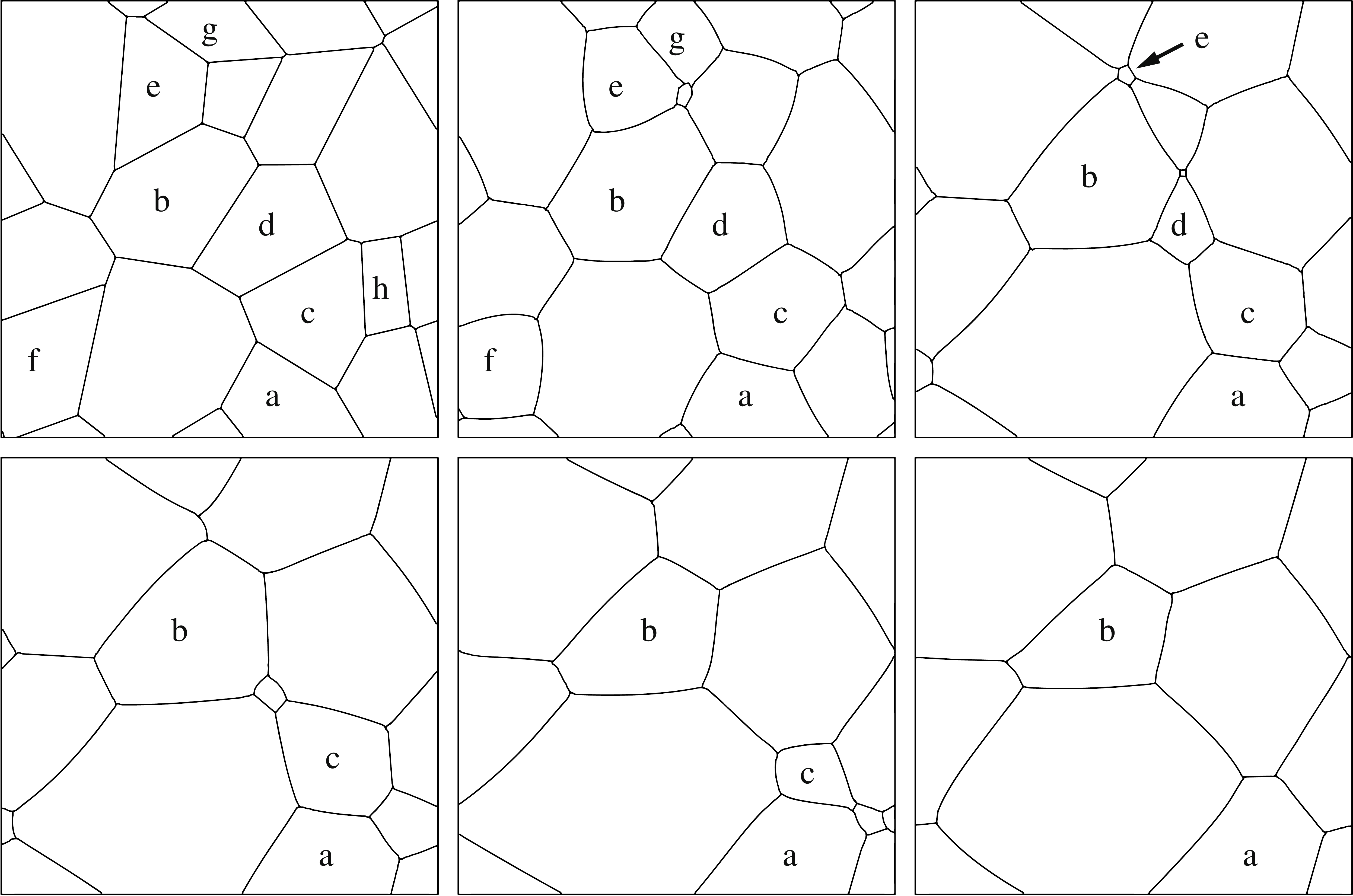}
\caption{Numerical result of the $15$-region system evolution under mean curvature.}
\label{15phases}
\end{center}
\end{figure}

\begin{figure}[p]
\begin{center}
\includegraphics[width=0.6\textwidth]{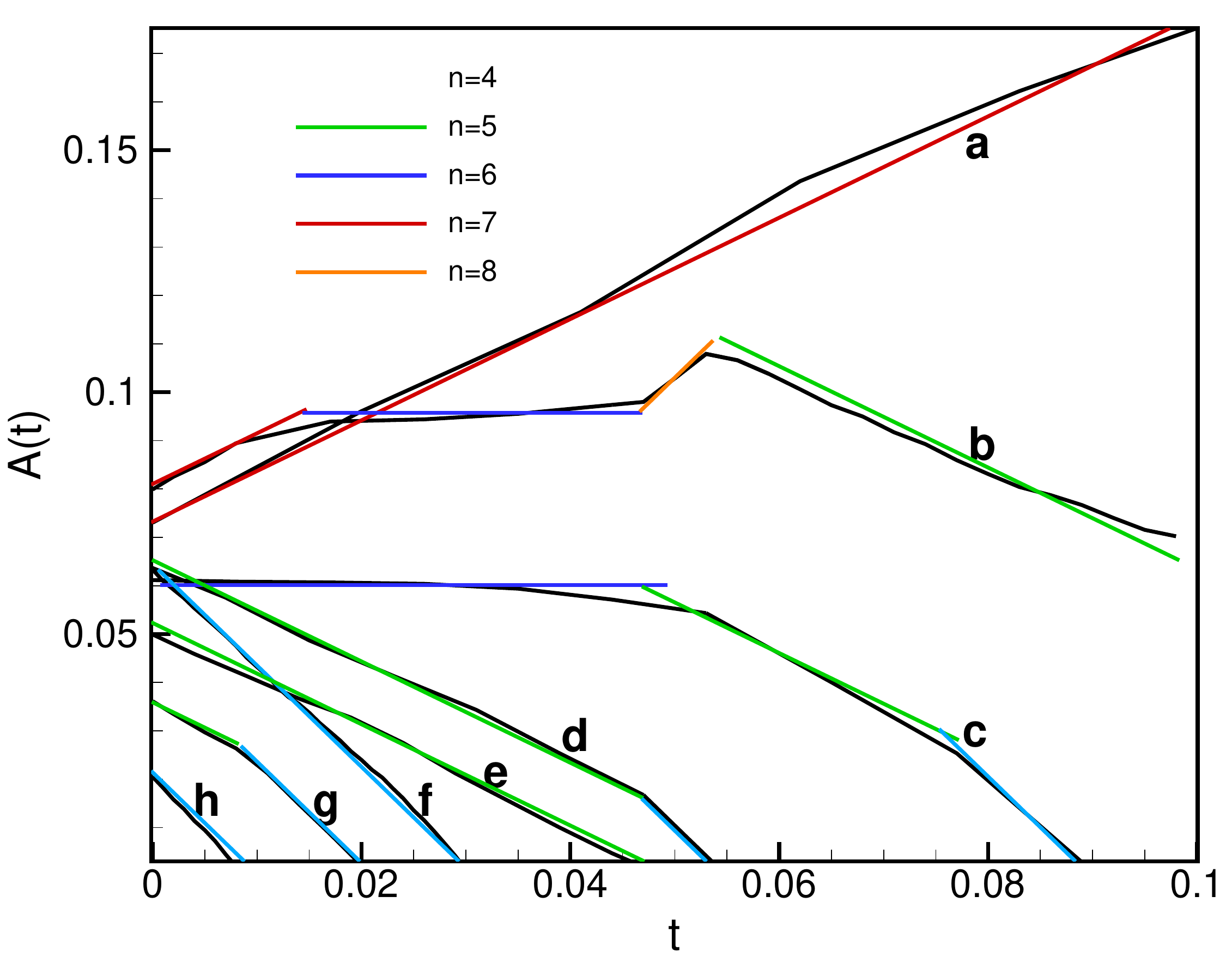}
\caption{Area as a function of time in the $15$-region system. The colored lines indicate von Neumann-Mullins' law in Eq. (\ref{Mullinslaw}).}
\label{area-t}
\end{center}
\end{figure}

We consider a case involving more regions and which serves to verify our method in a configuration with multiple triple points within a multi-region system. Initially $15$ regions are randomly placed in a domain $\left[0, 1\right]\times\left[0, 1\right]$ and evolve under a mean curvature generated velocity $\mathbf{u}=\kappa \mathbf{n}$. Periodic boundary conditions are imposed at the domain boundaries. The explicit Euler scheme is employed for temporal discretization. The time-step size is $\Delta t = \frac{h^2}{4}$ and the grid spacing is $h = \frac{1}{128}$. We want to verify that the von Neumann-Mullins' law is reproduced, which states that the rate of area $A$ growth or decay is a function of the number of edges $n$ of the phase \cite{saye2012analysis}
\begin{equation} \label{Mullinslaw}
\frac{dA}{dt}=2 \pi \gamma (\frac{n}{6}-1),
\end{equation}
where we set $\gamma = 1.0$. The initial number of regions $15$ is successively reduced to $6$ under the effect of mean curvature, as shown in Fig.\ref{15phases}. According to Fig.\ref{15phases},  region ``$h$'' initially has four edges and shrinks under mean curvature, leading to region destruction, consistent with von Neumann-Mullins' law. More specifically, the temporal growth and decay of the selected $8$ regions is shown in Fig. \ref{area-t} which exhibits a piecewise linear profile for each region, in agreement with von Neumann-Mullins' law which is indicated by a colored line in Fig. \ref{area-t}.

\subsection{Dry foam dynamics} \label{sec:foam}
As demonstration for a complex application we couple the high-resolution regional level-set method with the Navier-Stokes equation Eq. (\ref{ns_equ}) to simulate dry-foam dynamics.
\begin{equation} \label{ns_equ}
\frac{\partial \mathbf{U}}{\partial t}+\nabla \cdot \mathbf{F}_c=\nabla \cdot {{\mathbf{F}}_{v}}+\sigma \kappa \delta \mathbf{n},
\end{equation}
where $\mathbf{U} = (\rho$, $\rho u$, $\rho v$, $\rho w$, $\rho E)^T$, in which $\rho$, $u$, $v$, $\rho w$, and $\rho E$ are the density, the three velocity components, and the total energy, respectively. $\mathbf{F}_c$ and ${\mathbf{F}}_{v}$ are the convective and viscous flux tensor, respectively. The surface tension term $\sigma \kappa \delta \mathbf{n}$ on the right-hand side describes foam dynamics, with $\sigma$ being the surface tension and $\delta$ being a smoothed Dirac delta function. The surface tension force is calculated by Eq. (14) of Ref. \cite{saye2012analysis}. The Ohnesorge number is $\mathrm{Oh} = \mu/\sqrt{\rho \sigma L} = 0.00024$, i.e., the flow is surface-tension dominated. We note that the above surface-tension model may be imroved by Eq. (16) in Ref. \cite{saye2012analysis}, which is, however, not our main concern. Symmetry boundary conditions are applied. We use the 5th-order WENO scheme and a 4th-order central scheme for spatial discretization and a second-order strongly stable Runge-Kutta scheme for temporal discretization, both for solving N-S equation and level-set advection. The CFL number is $0.6$ and the grid size is $h=\frac{1}{256}$.

\begin{figure}[p]
\begin{center}
\includegraphics[width=0.9\textwidth]{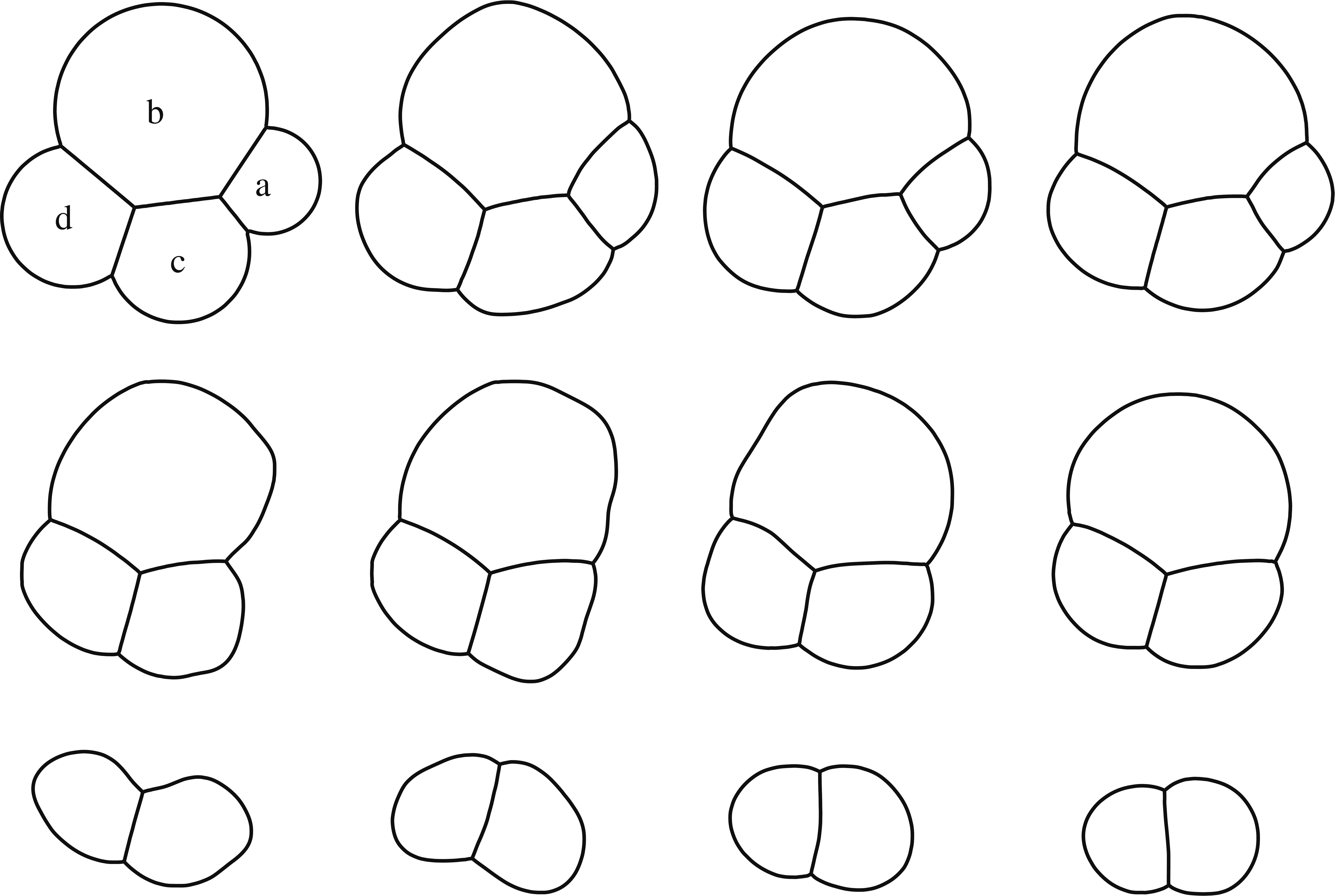}
\caption{Four-bubble cluster dynamics subject to sudden rupture. The breakup of regions $\Omega^{\chi_a}$ and $\Omega^{\chi_b}$ occurs at $t=3.0$ and $t=6.0$, respectively.}
\label{foam}
\end{center}
\end{figure}

In a dry-foam cluster the gas bubbles are separated by thin liquid films which corresponds to the interface network of our method. As we are not concerned with a specific physical problem, inertia and gravity effects of liquid membranes, gas exchange across permeable membranes, and Marangoni forces at the liquid-gas interface are neglected here for simplicity. We consider the interconnected membranes as massless and infinitely thin. Initially, four bubbles evolve from an artificial initial configuration to an equilibrium state, as shown in the first row of Fig.\ref{foam}. Following the breakup of $\Omega^{a}$, a second equilibrium state is attained under the surface tension force, see the second row of Fig.\ref{foam}. Then $\Omega^{b}$ breaks up, and a similar process is observed for the remaining two bubbles. Note that the bubble breakup is triggered explicitly, which is unphysical, it serves, however, our purpose of demonstrating that our method has the capability to capture the interface evolution of a foam cluster subjected to stimulated bubble breakup.

\section{Concluding remarks}\label{sec-conclusion}
The proposed method employs locally constructed level-set fields and the regional level-set method to overcome typical problems encountered with the numerical simulation of multi-region problems. 
The proposed numerical approach for multi-region problems and its algorithmic formulation have the following main properties: (1) As the proposed regional level-set method employs local signed level-set fields generated from the regional level-set field with a simple construction operator, it permits the implementation of high-resolution schemes for level-set transport in a straightforward way. (2) Instead of explicitly constructing the interface at every time step, we use a reconstruction operator to assemble the regional level set from multiple local level-set fields. This way we can ensure that the implicitly defined topology has no void or overlap artifacts. Moreover, a simple algorithm allows to distinguish different types of cells. It can be concluded from a range of test cases that the proposed method is more accurate than the Semi-Lagrangian regional level-set method. High order accuracy is demonstrated for some simple test cases where analytical results are known. Several increasingly complex configurations serve to demonstrate that improved accuracy and efficiency transfer to such test cases. A region-deconstruction example in foam dynamics demonstrates the feasibility of the method for complex applications.
\section*{Acknowledgment}
This work is partially supported by China Scholarship Council (No. 201306290030),
National Natural Science Foundation of China (No. 11628206) and 
Deutsche Forschungsgemeinschaft (HU 1527/6-1).
The project has received funding from the European Research Council (ERC) under the European Union's Horizon 2020 research and innovation program (No. 667483).
\bibliographystyle{model3-num-names}


\providecommand{\noopsort}[1]{}\providecommand{\singleletter}[1]{#1}%

\begin{algorithm}[p]
\caption{Cell type determination. Cell center data $\phi_{i,j}$ and $\chi_{i,j}$ are given, with indices $i=i_s,...,i_e, j=j_s,...,j_e, (i,j)\in \mathbb{N}$.}
\label{algo:type}
\begin{algorithmic}[1]
    \For{$i = i_s \ \mathbf{to} \ i_e$ and $j = j_s \ \mathbf{to} \ j_e$}
    \State for the cell $C_{i,j}$, the local index set for its square-shaped near neighborhood $V_{s}$ is $\mathrm{X}_s$ and the number of region indicators is $\mathcal{N}_s=n(\mathrm{X}_s)$
    \If {$\mathcal{N}_s = 1$}
       \State $C_{i,j}$ is a full cell because $1 \leq \mathcal{N}_C \leq \mathcal{N}_s =1$
    \ElsIf {$\mathcal{N}_s \geq 2$}
        \State construct the $\mathcal{N}_s$ local signed level-set fields in $V_{s}$:
            \State $\{\phi^r_{k,l} | \phi^r_{k,l} = \mathbf{C}_r(\varphi^{\chi}_{k,l}), r \in \mathrm{X}_s, C_{k,l} \in V_s\}$.
        \State Calculate the $\phi^r(\mathbf{x})$ at any points inside $C_{i,j}$ by bilinear interpolation interpolation.
        \State Check if $C_{i,j}$ is intersected by any region boundary $\partial \Omega^{\chi_r} = \{\mathbf{x} \in C_{i,j} | \phi^r(\mathbf{x}) = 0\}$ to determine the number of regions occupying the current cell, $\mathcal{N}_C =n(\mathrm{X}_{C_{i,j}})$:
            \State $\mathrm{X}_{C_{i,j}} = \left\{ { r \in \mathrm{X}_s | \partial \Omega^{\chi_r} \bigcap \partial \mathcal{C}_{i,j} \neq \varnothing}\right\}$
        \If {$\mathcal{N}_C = 1$}
            \State $C_{i,j}$ is a full cell
        \ElsIf {$\mathcal{N}_C = 2$}
            \State $C_{i,j}$ is a two-region cut cell
        \Else{}
            \State $C_{i,j}$ is a complex-region cut cell
        \EndIf
    \EndIf
    \EndFor
\end{algorithmic}
\end{algorithm}
\end{document}